\documentclass[a4paper]{amsart}

\usepackage{amsmath,amssymb}
\usepackage[dvipsnames]{xcolor}

\allowdisplaybreaks
\newtheorem{theorem}{Theorem}[section]
\newtheorem{proposition}[theorem]{Proposition}
\newtheorem{corollary}[theorem]{Corollary}

\newtheorem{question}[theorem]{Question}

\newtheorem{preremark}[theorem]{Remark}
\newtheorem{predefinition}[theorem]{Definition}
\newtheorem{preexample}[theorem]{Example}
\newtheorem{prenotation}[theorem]{Notation}
\newtheorem{preconjecture}[theorem]{Conjecture}
\newtheorem{assumption}[theorem]{Assumption}

\newenvironment{remark}{\begin{preremark}\rm}{\end{preremark}}
\newenvironment{definition}{\begin{predefinition}\rm}
	{\end{predefinition}}
\newenvironment{example}{\begin{preexample}\rm}{\end{preexample}}

\newcommand{\X}{{\mathbb{X}}}
\newcommand{\Y}{{\mathbb{Y}}}

\newcommand{\ZZ}{{\mathbb{Z}}}
\newcommand{\NN}{{\mathbb{N}}}

\newcommand{\Pbar}{{\overline{P}}}
\newcommand{\Rbar}{{\overline{R}}}

\let\epsilon=\varepsilon

\def\phi{{\varphi}}
\let\Psi=\varPsi
\let\Phi=\varPhi
\let\theta=\theta
\let\rho=\varrho

\def\HF{\mathop{\rm HF}\nolimits}
\def\HP{\mathop{\rm HP}\nolimits}

\def\Hom{\mathop{\rm Hom}\nolimits}
\def\Supp{\mathop{\rm Supp}\nolimits}

\newcommand{\charac}{\mathop{\rm char}\limits}

\newcommand{\Ker}{\mathop{\rm Ker}\nolimits}
\renewcommand{\Im}{\mathop{\rm Im}\nolimits}

\newcommand{\deh}{{\mathop{\rm deh}}}
\newcommand{\diff}{{\mathop{\rm diff}}}
\newcommand{\ri}{{\mathop{\rm ri}}}

\let\To=\longrightarrow

\def\tsum_#1^#2{{\textstyle\sum\limits_{#1}^{#2}}}

\begin{document}
\title[The Canonical Exact Sequence of Differential Modules for 0-Dim Schemes]
{The Canonical Exact Sequence of Differential Modules for 0-Dimensional Schemes}

\author{Tran N. K. Linh}
\address{Department of Mathematics,
	University of Education - Hue University, 34 Le Loi Street, Hue City, Vietnam}
\email{tnkhanhlinh@hueuni.edu.vn}

\author{Le Ngoc Long}
\address{Department of Mathematics,
	University of Education - Hue University, 34 Le Loi Street, Hue City, Vietnam}
\email{lelong@hueuni.edu.vn}

\date{\today}

\begin{abstract}
Given a 0-dimensional scheme $\X$ in $\mathbb{P}^n_K$ over a perfect field $K$, 
we examine the second differential power of its homogeneous vanishing ideal. 
This enables us to establish the canonical exact sequence 
for the associated K\"ahler differential module.
We also provide a formula for the Hilbert polynomial of  
K\"ahler differential modules when $\X$ is either a fat point scheme 
or a 0-dimensional locally monomial Gorenstein scheme.
\end{abstract}

\keywords{K\"ahler differential modules, differential powers,
	Hilbert functions, zero-dimensional schemes,  
	locally monomial Gorenstein schemes}

\subjclass[2020]{Primary 13N05; Secondary  13D40, 14N05, 13N10}

\maketitle

%
%
\section{Introduction}

Let $K$ be a perfect field, and
let $R$ be a finitely generated $K$-algebra of the form $R=P/I$,
where $P=K[X_0,\dots,X_n]$ is a polynomial ring over $K$ and 
$I$ is an ideal of $P$. One of powerful tools for investigating 
the $K$-algebra $R$ is its associated K\"ahler differential module.
This module is precisely defined as $\Omega^1_R = \Delta_R/\Delta_R^2$, where
$\Delta_R$ is the kernel of the multiplication map 
$\mu:R\otimes_KR\rightarrow R$ with $f\otimes g\mapsto fg$. 
The module $\Omega^1_R$ is finitely generated and it reflects 
many structural properties of the $K$-algebra $R$, including criteria for smoothness, 
regularity, and ramification (see, for instance, \cite{Gro1967, Kun1986, SS1974}). 

It is well-known that  $\Omega^1_R$ admits
the following canonical exact sequence of $R$-modules 
$$
I/I^2 \stackrel{\delta}{\longrightarrow} \Omega^1_P/ I\Omega^1_P
\stackrel{\gamma}{\longrightarrow} \Omega^1_R \longrightarrow 0
$$ 
where $\delta(f+I^2) = df + I\Omega^1_P$ for each $f\in I$ and 
the map $\gamma$ is induced from the functorial map $\Omega^1_P \rightarrow \Omega^1_R$.
A significant result by Mohan Kumar \cite{Kum1978} demonstrates that 
when $R$ is reduced, $I$ is a complete intersection if and only if 
the exact canonical sequence forms a free resolution of $\Omega^1_R$, 
i.e., if and only if 
$$
0 \longrightarrow I/I^2 \stackrel{\delta}{\longrightarrow} 
R^{n+1} \stackrel{\gamma}{\longrightarrow} 
\Omega^1_R \longrightarrow 0
$$ 
is exact and $I/I^2$ is a free $R$-module. 
This exactness condition has been further explored in subsequent works.
In \cite{DK1999}, De Dominicis and Kreuzer proved that 
if $I = I_\X$ is the homogeneous vanishing ideal of a finite set of points 
$\X = \{p_1, \dots, p_s\}$ 
in $\mathbb{P}^n_K$ over a field $K$ of characteristic zero, 
then the sequence of graded $R$-modules
\begin{equation}\label{CanonicalExactSeq}
0 \longrightarrow I_\X/I_{\Y} \stackrel{\delta}{\longrightarrow} R^{n+1}(-1) 
\stackrel{\gamma}{\longrightarrow} \Omega^1_R \longrightarrow 0
\end{equation}
is exact, where $I_\Y$ is the homogeneous vanishing ideal of 
the double point scheme $\Y$ in $\mathbb{P}^n_K$. 
Recall that, for $K$-rational points $p_1,\dots,p_s$ and positive integers
$m_1,\dots, m_s$, the saturated homogeneous ideal 
$I_\X = I_{p_1}^{m_1}\cap\cdots\cap I_{p_s}^{m_s}$
of $P$ defines a fat point scheme $\X= m_1p_1+\cdots+m_sp_s$ in $\mathbb{P}^n_K$.
Recently, when $\charac(K)=0$, 
the canonical exact sequence \eqref{CanonicalExactSeq}
was generalized in \cite[Theorem 1.7]{KLL2015} to such a fat point scheme
$\X= m_1p_1+\cdots+m_sp_s$ using its fattening $\Y=(m_1+1)p_1+\cdots+(m_s+1)p_s$ 
in~$\mathbb{P}^n_K$.
Inspired by these advancements, it is natural to ask the following question:

\begin{question}\label{Question 1}
	Given an arbitrary 0-dimensional scheme $\X$ in $\mathbb{P}^n_K$
	over a perfect field $K$ of characteristic $p\ge 0$, 
	does a subscheme $\Y$ of $\mathbb{P}^n_K$ exist for which 
	the canonical exact sequence \eqref{CanonicalExactSeq} for the K\"ahler differential 
	module for $\X$ remains valid?
\end{question}

To answer this question, we first deal with the differential powers 
of the homogeneous vanishing ideal $I_\X$ of the scheme $\X$. 
For $m\ge 0$, let $D_P^m$ be the $P$-module of $K$-linear 
differential operators on $P$ of order at most~$m$. 
On account of \cite[Definition 2.2]{DDGHN}, the $m$-th differential power 
of $I_\X$ is given by $I_\X^{\langle m\rangle} = \{f\in P\mid \delta(f)\in I_\X\ 
\text{for all $\delta\in D^{m-1}_P$}\}$.
For our purpose, we focus on investigating the second differential power
$I_\X^\diff := I_\X^{\langle 2\rangle}$.
Our first result is the following theorem that shows how to compute
a finite generating system of $I_\X^\diff$, when $\X$ contains only a
monomial point. 
\begin{theorem}[Theorem \ref{Prop-Diff-MonomialIdeal}]\label{Theorem 1.2}
	Let $K$ be a perfect field of characteristic $p\ge 0$. 
	Assume that the homogeneous vanishing ideal $I_\X$ 
	of the scheme $\X$ is a proper monomial ideal of $P$, minimally generated by 
	terms $t_1,\dots,t_r \in \mathbb{T}^{n+1}$. 
	\begin{enumerate}
		\item[(a)] The ideal $I_\X^\diff$ is a monomial subideal of $I_\X$.
		
		\item[(b)] For $1\le j\le r$, we write 
		$t_j = X_0^{\alpha_{0j}}\cdots X_n^{\alpha_{nj}}$, 
		and for $2\le k\le n+1$ and $1\le j_1\le\cdots\le j_k\le r$, we define
		$t_{j_1,\dots,j_k} := X_0^{\beta_0}\cdots X_n^{\beta_n}$, where, for $i=0,\dots,n$,
		\begin{equation} 
			\beta_i \;:=\; 
			\begin{cases}
				0 & \text{if $\alpha_{ij_1}=\dots= \alpha_{ij_k}=0$},\\
				\alpha_{ij_1} & 
				\text{if $\alpha_{ij_1}=\dots= \alpha_{ij_k}>0$ and $p\mid \alpha_{ij_1}$},\\
				\alpha_{ij_1}+1 & 
				\text{if $\alpha_{ij_1}=\dots= \alpha_{ij_k}>0$ and $p\nmid \alpha_{ij_1}$},\\
				\max\{\alpha_{ij_1},\dots,\alpha_{ij_k}\} & \text{otherwise}.
			\end{cases}
		\end{equation}
		Then 
		$$
		I_\X^\diff \;=\; \langle\, t_{j_1,\dots,j_k} \mid 
		1\le j_1\le\cdots\le j_k\le r;\; 2\le k\le n+1 \,\rangle.
		$$
	\end{enumerate}
\end{theorem}
In the general case of a 0-dimensional scheme $\X$, we show that 
$I_\X^\diff$ is indeed a saturated homogeneous ideal in $P$, 
and so it defines a 0-dimensional scheme $\Y$ in~$\mathbb{P}^n_K$
which is called the \textit{differential fattening} of~$\X$. 
By $x_i$ we denote  the image of $X_i$ in $R=P/I_\X$ for $i=0,\dots,n$.
Under this terminology,
we give an affirmative answer for Question~\ref{Question 1}.
\begin{theorem}[Theorem \ref{Thm-CanonicalExactSequence}]\label{Theorem 1.3}
	Let $\X$ be a 0-dimensional scheme in $\mathbb{P}^n_K$ with homogeneous 
	vanishing ideal~$I_\X$, and let $\Y$ be the differential fattening of $\X$
	with its homogeneous vanishing ideal~$I_\Y$.
	The sequence of graded $R$-modules
	$$
	0 \longrightarrow I_\X/ I_\Y \stackrel{\delta}{\longrightarrow}
	R^{n+1}(-1)  \stackrel{\gamma}{\longrightarrow} \Omega^1_R
	\longrightarrow 0
	$$
	is exact, where $\delta(f + I_\Y) 
	=\sum_{i=0}^{n}\frac{\partial f}{\partial x_i}e_i$ 
	for every homogeneous polynomial $f\in I_\X$ and
	$\gamma(e_i) = dx_i$ for $i=0,\dots,n$.
\end{theorem}
Clearly, the theorem enables us to compute 
the Hilbert function of the K\"ahler differential module $\Omega^1_R$
and to bound its regularity index via those of $\X$ and of~$\Y$.
When $\X=m_1p_1+\cdots+m_sp_s$ is a fat point scheme in $\mathbb{P}^n_K$,
\cite[Remark 7.1]{KLL2025} poses the question: \textit{does the Hilbert polynomial 
of the $m$-form K\"ahler differential module $\Omega^m_R$
depend only on $m$, $n$ and the multiplicities $m_1,\dots,m_s$?}
A positive answer for this question is given by \cite[Theorem 7.12]{KLL2025}
when either $\charac(K)=0$ or $\charac(K)>\max\{m_j\mid j=1,\dots,s\}$.
Based on Theorems~\ref{Theorem 1.2} and~\ref{Theorem 1.3},
we give a formula for the Hilbert polynomial of $\Omega^1_R$ 
and this improves \cite[Theorem 7.12]{KLL2025} for $m=1$ and 
$0\le \charac(K) \le \max\{m_j\mid j=1,\dots,s\}$ (see Corollary~\ref{Cor-HFOmega^1_R}). 
Unfortunately, we do not know an answer to this question 
when $m\ge 2$ and $0\le \charac(K) \le \max\{m_j\mid j=1,\dots,s\}$.

On the other hand, we are also interested in finding a formula for
the Hilbert polynomial of the $m$-form K\"ahler differential module 
for 0-dimensional locally monomial Gorenstein schemes in $\mathbb{P}^n_K$.
Here, a 0-dimensional scheme $\X$ with support $\Supp(\X)=\{p_1,\dots,p_s\}$ 
in $\mathbb{P}^n_K$ is a locally monomial Gorenstein scheme if, for each 
$j\in\{1,\dots,s\}$, the local ring $\mathcal{O}_j$ of $\X$ at $p_j$ satisfies
$\mathcal{O}_j \cong K[X_1,\dots,X_n]/\mathfrak{Q}_j$, where $\mathfrak{Q}_j$
is a monomial Gorenstein ideal generated by pure powers of the indeterminates. 
The class of these schemes includes the class 
of 0-dimensional curvilinear schemes which have been shown to play 
an important role in the proof of the so-called Alexander-Hirschowitz theorem 
(see, for instance, \cite{Chan2000}). We derive the following corollary.
\begin{corollary}[Corollary \ref{Cor-HPofOmage^m_R-LMG-Sch}] \label{Corollary 3}
	Let $\X$ be a 0-dimensional locally monomial Gorenstein scheme in $\mathbb{P}^n_K$
	with support $\Supp(\X)=\{p_1,\dots,p_s\}$, and let $p=\mathrm{char}(K)$.
	For $j=1,\dots s$, we write $\mathcal{O}_j \cong K[X_1,\dots,X_n]/\mathfrak{Q}_j$,
	where $\mathfrak{Q}_j = \langle X_1^{k_{1j}},\dots,X_n^{k_{nj}} \rangle$ 
	with $k_{ij}\ge 1$. For $1\le i_1<\cdots<i_m\le n$ and $1\le j\le s$, let 
	$$
	\Gamma_{(i_1,\dots,i_m),j} = \{i\in \{i_1,\dots,i_m\} :\;\,	p \nmid k_{ij}\}.
	$$
	For $m\ge 1$, we have 
	\begin{align*}
		\HP(\Omega^m_R) &= \sum_{j=1}^s
		\Big[ \sum_{1\le i_1<\dots<i_{m-1}\le n} 
		\prod_{i\notin \Gamma_{(i_1,\dots,i_{m-1}),j}}k_{ij} 
		\prod_{i\in\Gamma_{(i_1,\dots,i_{m-1}),j}} (k_{ij}-1) \\
		&\quad +\sum_{1\le l_1<\dots<l_m\le n} 
		\prod_{i\notin \Gamma_{(l_1,\dots,l_m),j}}k_{ij}
		\prod_{i\in\Gamma_{(l_1,\dots,l_m),j}} (k_{ij}-1)\Big].
	\end{align*}
\end{corollary}

This paper is outlined as follows. In Section 2, we first recall the needed facts about
the K\"ahler differential modules, the ring of $K$-linear differential 
operators, and the $K$-linear differential powers.  Then we consider these invariants
for a 0-dimensional scheme in $\mathbb{P}^n_K$. 
In Section 3, we prove Theorems~\ref{Theorem 1.2} and~\ref{Theorem 1.3} and
derive their consequences (see Corollaries \ref{Cor-Diff-MonomialIdeal}, 
\ref{Cor-HFOmega^1_R}, and \ref{Cor-Curvilinear}). In the final section, we establish
a formula for the Hilbert polynomial of the $m$-form K\"ahler differential module
associated to a 0-dimensional locally monomial Gorenstein scheme 
(see Proposition~\ref{Prop-HPOfOmega^1_R-LMG-Scheme} and Corollary~\ref{Corollary 3}).

Unless explicitly stated otherwise, we adhere to the definitions and notation
introduced in \cite{KR2000, KR2005} and \cite{Kun1986}. 
All examples in this paper were calculated by using the
computer algebraic system \texttt{ApCoCoA} \cite{ApCoCoA}.

\bigbreak
%
%

\section{Differential Modules and Differential Powers of Ideals}
\label{Section 2}

In the following, let $K$ be a perfect field and let $P=K[X_0,\dots,X_n]$ 
be a polynomial ring in indeterminates $X_0,\dots,X_n$ over $K$.
Let $R$ be a finitely generated $K$-algebra of the form $R = P/I$,
where $I$ is an ideal in $P$. 
The image of $X_i$ in $R$ is denoted by $x_i$ for $i=0,\dots,n$.  
We denote by $\mu$ the multiplication map $\mu: R \otimes_K R \To R$ 
given by $\mu(f\otimes g)=fg$ for $f,g\in R$. 
This is a $K$-algebra homomorphism and its kernel $\Delta_R :=\Ker(\mu)$
is an ideal of the algebra $R \otimes_K R$. In particular, 
the ideal $\Delta_R$ is generated by 
$\{ 1\otimes x_i - x_i\otimes 1 \mid i\in\{0,\dots,n\}\}$.

\begin{definition}\label{def:KDdefn}
	Let $R = P/I$ be a finitely generated $K$-algebra.
	\begin{enumerate}
		\item[(a)] The finitely generated $R$-module $\Omega^1_{R/K}=\Delta_R /\Delta_R^2$ is called the {\bf K\"ahler differential module}
		(or the {\bf module of K\"ahler differentials}) of~$R/K$.
		
		\item[(b)] For every $m\ge 0$, the $m$-th exterior power 
		$\Omega^m_{R/K} := \Lambda^m_R \Omega^1_{R/K}$ 
		is called the {\bf $m$-form K\"ahler differential module}  
		(or the {\bf module of K\"ahler differential $m$-forms}) of~$R/K$.
		
		\item[(c)] The exterior algebra $\Omega^{\bullet}_{R/K} :=
		\bigoplus_{m\ge 0} \Omega^m_{R/K}$ is called the {\bf K\"ahler differential algebra} 
		of~$R/K$.
	\end{enumerate}	
\end{definition}

Unless explicitly stated otherwise, all algebras will be $K$-algebras.
So, we usually simplify the notation and write $\Omega^m_{R}$ instead of
$\Omega^m_{R/K}$ for all $m\ge 0$, as well as $\Omega^{\bullet}_{R}$ instead of 
$\Omega^{\bullet}_{R/K}$. Notice that $\Omega^0_{R} = R$ and if $I$ is a
homogeneous ideal of $P$, then $\Omega^m_{R}$ is a graded $R$-module 
for every $m\ge 0$. The module of K\"ahler differentials of $R/K$ comes
together with the {\it universal derivation} $d_{R}: R \rightarrow \Omega^1_R$ 
which is defined by $d_{R}(f) = f\otimes 1 - 1 \otimes f + \Delta_R^2$
for $f\in R$. We simply write $d$ for $d_{R}$, if no confusion can arise. 

As in \cite[Section 1]{Kun1986}, the set of all $K$-linear derivations 
$\partial: R\rightarrow R$ is an $R$-module and is denoted by $\operatorname{Der}_K(R)$.
By \cite[Proposition 1.23]{Kun1986}, there is an isomorphism 
$$
\Hom_R(\Omega^1_{R},R) \To \operatorname{Der}_K(R),\; 
\varphi \longmapsto \varphi\circ d.
$$ 
In addition, both $\operatorname{Der}_K(R)$ and $\Hom_{R}(R,R)$
are submodules of $\Hom_{K}(R,R)$. The latter, $\Hom_{K}(R,R)$, 
is the ring of endomorphisms of $R/K$, with operations of
sum and composition of endomorphisms.
If $\delta, \theta \in \Hom_{K}(R,R)$, then their commutator is defined as
$$
[\delta,\theta] := \delta\circ \theta - \theta\circ\delta.
$$
Obviously, $[\delta,\theta]\in \Hom_{K}(R,R)$. Every element in $\Hom_{R}(R,R)$ 
is simply a multiplication map $\mu_f: R\longrightarrow R$
for some $f\in R$, where $\mu_f(g)=fg$ for all $g\in R$. 
We write $[\delta, f]$ for $[\delta, \mu_f]$, and so that
$$
[\delta, f] = \delta\circ \mu_f - \mu_f\circ\delta
\quad\text{and}\quad  
[f,g] = 0,
$$
for all $f,g\in R$ and $\delta\in \Hom_{K}(R,R)$.

Based on \cite[Section 2]{HS1969} or \cite[Chapter 3, Section 1]{Cout1995}, 
we have the following notion.

\begin{definition}
Let $R=P/I$ be a finitely generated $K$-algebra, and let $m\in \NN$.
The \textbf{$K$-linear differential operators of $R$ of order at most $m$},
$D^m_R\subseteq \Hom_{K}(R,R)$, are defined inductively as follows:
\begin{enumerate}
	\item[(i)] $D^0_R = \Hom_{R}(R,R)$;
			
	\item[(ii)] $D^m_R = \{\delta \in \Hom_{K}(R,R) 
	\mid [\delta, f] \in D^{m-1}_R 	\ \text{for every $f\in R$} \}$.
\end{enumerate}
The \textbf{ring of $K$-linear differential operators} of $R$ is defined by 
$
D_R = \bigcup_{m\in\NN} D^m_R.
$
\end{definition} 

Note that the ring structure on $D_R$ is also given by composition and satisfies
$D^m_K\subseteq D^{m+1}_R$ for $m\ge 0$ and $D^m_R\circ D^k_R \subseteq D^{m+k}_R$ 
(see, e.g., \cite[Chapter 3, Proposition 1.2]{Cout1995}).
The ring $D_R$ is a non-commutative ring, 
and it is in general not Noetherian (see \cite{BGG1972}).

When $R=P=K[X_0,\dots,X_n]$ and $K$ is of characteristic zero, 
the ring $D_P$ is known as the Weyl algebra of $P$ and 
it has a concrete description below (see, for instance,
\cite[Chapter I, Proposition 1.2]{Bjo1979} and \cite[Theorem 2]{Hart1983}). 
We refer \cite{Smith1986} for further information of $D_P$ in positive characteristic. 
In the following, $\partial_0,\dots,\partial_n$ are
the partial derivatives defined by $\partial_i(f)=\frac{\partial f}{\partial X_i}$ for $f\in P$. 
For $\alpha=(\alpha_0,\dots,\alpha_n)\in \NN^{n+1}$, by $X^\alpha$ 
we mean the term $X_0^{\alpha_0}\cdots X_n^{\alpha_n}$ and its degree is 
$|\alpha|= \alpha_0+\cdots+\alpha_n$, and similarly,
for $\beta=(\beta_0,\dots,\beta_n)\in \NN^{n+1}$, by $\partial^\beta$ 
we denote a $\partial$-term $\partial_0^{\beta_0}\cdots \partial_n^{\beta_n}$.

\begin{proposition}\label{Prop-TheWeylAlgeOfPolyRing}
If $K$ is a field of characteristic zero, then 
$D_P = P[\operatorname{Der}_K(P)] = K[X_0,\dots,X_n,\partial_0,\dots,\partial_n]$ 
and the set $\{X^\alpha \partial^\beta \mid \alpha,\beta \in \NN^{n+1} \}$ 
is a $K$-basis of $D_P$.
In particular, any $\delta \in D_P^m$ can be uniquely written as
$$
\delta = \sum_{|\alpha|\le m} f_{\alpha} \partial^{\alpha}\qquad (f_\alpha \in P).
$$
\end{proposition}

For the finitely generated $K$-algebra $R=P/I$, there is an $(R\otimes_KR)$-module 
structure on $\Hom_{K}(R,R)$ given by the rule
$$
((f\otimes g)\cdot \theta)(h) = f\theta(gh)
$$
for $f,g,h\in R$ and $\theta\in \Hom_{K}(R,R)$.
Via this module structure, we have the following characterization 
of differential operators (see \cite[Lemma 2.2.1 and Proposition 2.2.3]{HS1969}).

\begin{proposition}
The set $D^m_R$ of $K$-linear differential operators of order $\le m$  
is an $(R\otimes_K R)$-submodule (and so an $R$-submodule) of $\Hom_{K}(R,R)$. 
Moreover, $\delta\in D_R^m$ 
if and only if $\Delta_R\cdot\delta \subseteq D_R^{m-1}$
if and only if $\Delta_R^{m+1}\cdot\delta =0$.
\end{proposition}

The following notion of differential powers of an ideal in $R$
was recently introduced in \cite[Definition 2.2]{DDGHN}.

\begin{definition}
Let $J$ be an ideal of $R$ and let $m$ be a positive integer.
The \textbf{$m$-th $K$-linear differential power} of $J$ is defined by
$$
J^{\langle m\rangle} := \{\, f\in R \mid \delta(f) \in J \ 
\text{for all $\delta\in D^{m-1}_{R}$} \,\}.
$$
\end{definition}

It is straightforward to verify that $J^{\langle m\rangle}$ is a subideal
of $J$ in~$R$ and if $J_1\subseteq J_2$ then 
$J_1^{\langle m\rangle}\subseteq J_2^{\langle m\rangle}$.
The next proposition gathers some additional properties 
of differential powers from \cite[Section 2]{DDGHN}.

\begin{proposition}\label{Prop-DifferentialPowers}
Let $J, J_1, J_2$ be ideals of the finitely generated $K$-algebra $R$, 
and let $m$ be a positive integer. Then 
\begin{enumerate}
	\item[(a)] $J^m\subseteq J^{\langle m\rangle} \subseteq J$;
		
	\item[(b)] $(J_1\cap J_2)^{\langle m\rangle} = J_1^{\langle m\rangle}\cap J_2^{\langle m\rangle}$;
		
	\item[(c)] If $R=P$ and $J$ is a maximal ideal of $P$, then 
		$J^{\langle m\rangle} = J^m$.	
\end{enumerate}
\end{proposition}

In what follows, let $\X$ be a 0-dimensional scheme in the $n$-dimensional projective 
space~$\mathbb{P}^n_K$ over the perfect field $K$, and let $I_\X$ be
the saturated homogeneous vanishing ideal of $\X$ in $P=K[X_0,\dots,X_n]$.
Then the homogeneous coordinate ring of $\X$ is $R= P/I_\X$. 
The ring~$R$ is a 1-dimensional, standard graded, Cohen-Macaulay $K$-algebra. 
The set of closed points of $\X$ is known as the support of $\X$ and 
is denoted by $\Supp(\X)=\{p_1,\dots,p_s\}$.

\begin{assumption}\label{Assumption 2.7}
No point of the support of~$\X$ lies on the hyperplane $D_+(X_0)$.
\end{assumption}

Under this assumption, the elements $x_0$ and $x_0-1$ are non-zerodivisors of $R$, 
the ring $\Rbar = R/ \langle x_0\rangle$ is a 0-dimensional local graded $K$-algebra, 
especially, $\Rbar$ is a finite dimensional $K$-vector space.
Moreover, $\X$ is contained in the affine space $\mathbb{A}^n \cong D_+(X_0)$.
The affine coordinate ring of~$\X$, viewed as a subscheme of~$\mathbb{A}^n$, is given by
$S = R/ \langle x_0- 1\rangle \cong A/J_\X$, where $A=K[X_1,\dots,X_n]$
and $J_\X=I_\X^\deh$ is the dehomogenization of~$I_\X$ with respect to~$X_0$. 
The ring~$S$ is a 0-dimensional affine $K$-algebra 
and hence a finite dimensional $K$-vector space. 

\begin{corollary}\label{Cor-Diff-Powers}
Let $\X$ be a 0-dimensional scheme in $\mathbb{P}^n_K$ 
with support $\Supp(\X)=\{p_1,\dots,p_s\}$,
and let $I_{p_j}$ be the homogeneous vanishing 
ideal of~$p_j$ for $j\in\{1,\dots,s\}$\!.
For $m\ge 1$, we have 
$$
I_\X^{\langle m\rangle} = I_{p_1}^{\langle m\rangle}\cap\cdots\cap I_{p_s}^{\langle m\rangle}
$$
and $\sqrt{I_{p_j}^{\langle m\rangle}} = \sqrt{I_{p_j}}$ for $j\in\{1,\dots,s\}$.
\end{corollary}

\begin{proof}
This follows from the primary decomposition $I_\X= I_{p_1}\cap\cdots\cap I_{p_s}$ 
and by Proposition~\ref{Prop-DifferentialPowers}(a)-(b).
\end{proof}

The module of K\"ahler differential $m$-forms for the scheme $\X$ is 
the finitely generated graded $R$-module $\Omega^m_{R}$ and 
the K\"ahler differential algebra for the scheme $\X$ is the exterior algebra 
$\Omega^{\bullet}_{R}=\bigoplus_{m\ge 0} \Omega^m_{R}$.

\begin{proposition}\cite[Proposition 4.2]{KLL2025} \label{prop:PresentationOmegaRm}
Let $\X$ be a 0-dimensional scheme in $\mathbb{P}^n_K$ with homogeneous
coordinate ring $R=P/I_\X$.
\begin{enumerate}
\item[(a)] If $m\ge n+2$, then $\Omega^m_R = \langle 0 \rangle$.
		
\item[(b)] For $m=0,\dots,n+1$, the $P$-module $\Omega^m_P$ is a free
$P$-module of rank $\binom{n+1}{m}$ with basis
$\{ dX_{i_1} \wedge \cdots \wedge
dX_{i_m} \mid 0\le i_1 < \cdots < i_m \le n \}$.
		
\item[(c)] For $m\ge 1$, the module of K\"ahler differential $m$-forms satisfies
$$
\Omega^m_R \;\cong\;   \Omega^m_P\;/ \big(  I_\X \Omega^m_P +
dI_\X \wedge \Omega^{m-1}_P  \big)
$$
where $dI_\X$ is the $P$-submodule of~$\Omega^1_P$ generated by
$\{df \mid f\in I_\X\}$.

\item[(d)] Let $I_\X = \langle f_1,\dots,f_k \rangle$ for some $f_1,\dots,f_k\in P$.
Then 
$$
d I_\X = \langle df_1,\dots,df_k \rangle + \langle f_i dX_j
\mid i\in \{1,\dots,k\},\; j\in \{0,\dots,n\} \rangle.
$$
\end{enumerate}
\end{proposition}

\begin{definition}
Let $M=\bigoplus_{i\in\ZZ} M_i$ be a finitely generated graded $R$-module.
\begin{enumerate}
	\item[(a)] The map $\HF_M:\; \ZZ \To \ZZ$ given by $\HF_M(i)= \dim_K(M_i)$
	for all $i\in\ZZ$ is called the {\bf Hilbert function} of~$M$. 
	The \textbf{Hilbert polynomial} of $M$ exists and is denoted by $\HP(M)$.
	
	\item[(b)] The number 
	$\ri(M) = \min\{i\in\ZZ \mid \HF_M(j) = \HP_M(j), \forall j\ge i\}$
	is called the {\bf regularity index} of~$M$.
\end{enumerate}
\end{definition}

\begin{remark}
The function $\HF_\X(i):=\HF_R(i)$ for all $i\in\ZZ$ is called
the \textbf{Hilbert function of~$\X$} and the number $r_\X := \ri(R)$ is called
the \textbf{regularity index of $\X$}. We have
$$
1 = \HF_\X(0) < \HF_\X(1) < \cdots < \HF_\X(r_\X) = \HF_\X(r_\X + 1) = \cdots = \deg(\X),
$$
where $\deg(\X)$ is the degree of $\X$, i.e., the constant Hilbert polynomial of $R$.
\end{remark}

It is worth noting that $\deg(\X)=\dim_K(S)$.
For $i\ge r_\X$, the map $\epsilon_i:\; R_i \To S$ defined by $f\mapsto f^\deh$ 
is an isomorphism of $K$-vector spaces (see \cite{KLL2019}). This structure map
induces a valuable connection between the modules $\Omega^m_R$ and $\Omega^m_{S}$
(see \cite{KLL2025}), ultimately leading to the proof of the following result.

\begin{proposition}\cite[Proposition 5.4]{KLL2025}\label{prop:HigherKDM}
Let $\X$ be a 0-dimensional scheme in~$\mathbb{P}^n$ with
homogeneous coordinate ring $R$ and affine coordinate ring 
$S=R/\langle x_0-1\rangle$,	and let $m\ge 1$. Then
$$
\HP(\Omega^m_R) = \dim_K(\Omega^m_{S}) + \dim_K(\Omega^{m-1}_{S})
\quad \text{and}\quad 
\ri(\Omega^m_R) \le 2 r_\X + m.
$$
In addition, if $\charac(K) = 0$ or if $\charac(K)=p>0$ and $p \nmid (2r_\X+n+1)$, 
then $\ri(\Omega^{n+1}_{R})\le 2 r_\X+n$.
\end{proposition}

\bigbreak
%
%

\section{The Canonical Exact Sequence for K\"ahler Differentials}
\label{Section 3}

As in the previous section, let $K$ be a perfect field, 
let $\X$ be a 0-dimensional scheme in $\mathbb{P}^n_K$
with homogeneous vanishing ideal $I_\X$ and homogeneous coordinate 
ring $R=P/I_\X$. Further, let $\Supp(\X)=\{p_1,\dots,p_s\}$, and let 
$I_{p_j}$ be the homogeneous vanishing ideal of the point $p_j$ for every 
$j\in\{1,\dots,s\}$. 
Our goal in this section is to give an answer to Question~\ref{Question 1}
posted in Section 1, regarding the existence of the canonical exact sequence
for the K\"ahler differential module $\Omega^1_R$ of the scheme $\X$.
To achieve this goal, we will first examine the second differential powers of
the homogeneous vanishing ideal $I_\X$ of the scheme $\X$. 

Let $\mathfrak{M} = \langle X_0,\dots,X_n\rangle\subseteq P$ 
and let $D_P$ denote the ring of $K$-linear differential operators on $P$. 
The second differential power $I_\X^{\langle 2\rangle}$ of $I_\X$ is given by 
$$
I_\X^{\langle 2\rangle} \; =\; \{\, f \in P \mid \delta(f) \in I_\X\
\text{for all $\delta\in D^1_P$} \,\}.
$$
This ideal is also written as $I_\X^{\diff}$. For $i=0,\dots,n$, we will
use $\partial_i$ to represent~$\frac{\partial }{\partial X_i}$.

\begin{proposition}\label{Prop-Diff-HomogVanishingIdeal}
Let $\X$ be a 0-dimensional scheme in $\mathbb{P}^n_K$ 
with support $\operatorname{Supp}(\X)=\{p_1,\dots,p_s\}$,
and let $I_{p_j}$ be the homogeneous vanishing 
ideal of $p_j$ for $j\in\{1,\dots,s\}$\!. 
\begin{enumerate}
\item[(a)] The ideal $I_\X^\diff$ is a saturated homogeneous ideal 
of~$P$ and 
\begin{equation}\label{Equ-DiffIdeal}
	I_\X^\diff \; =\; \{\, f \in I_\X \mid \partial_i f \in I_\X\
	\text{for all $i\in\{0,\dots,n\}$} \,\} \subseteq P.
\end{equation}
		
\item[(b)] We have $I_\X^\diff = I_{p_1}^\diff\cap\cdots
\cap I_{p_s}^\diff$, and $\sqrt{I_{p_j}^\diff}=\sqrt{I_{p_j}}$
for $j=1,\dots,s$.
\end{enumerate}
\end{proposition}

\begin{proof}
(a)\quad According to \cite[Chapter 3, Proposition 1.1]{Cout1995}, 
we first observe that $D^1_P = P + \operatorname{Der}_K(P)$, 
where the elements of $P$ are identified with their multiplication operators 
of order zero.  By \cite[Example 1.22(c)]{Kun1986}, $\operatorname{Der}_K(P)$
is generated as an $R$-module by the partial derivatives $\partial_0,\dots,\partial_n$.
They form in fact a basis of $\operatorname{Der}_K(P)$ since 
if $\sum_{i=0}^{n}f_i \partial_i =0$ for $f_0,\dots,f_n\in P$ then
$f_k = \sum_{i=0}^{n}f_i \partial_i(X_k) =0$ for $k=0,\dots,n$.
Thus we have  $D^1_P = P\oplus P\partial_0 \oplus\cdots\oplus P\partial_n$.
Consequently, the equality \eqref{Equ-DiffIdeal} holds true for $I_\X^\diff$.

Now let $g = g_0+\cdots+g_d\in I_\X^\diff$ be a nonzero polynomial 
of degree $d$, where $g_k\in P_k$ is the $k$-th homogeneous component of $g$.
For each $i\in\{0,\dots,n\}$, we have 
$\partial_i(g) = \partial_i(g_0)+\cdots+\partial_i(g_d)\in I_\X$. Since $I_\X$ is 
homogeneous, $\partial_i(g_j) \in I_\X$ for $j\in\{0,\dots,d\}$.
Thus $g_j\in I_\X^\diff$ for $j\in\{0,\dots,d\}$, 
and hence $I_\X^\diff$ is homogeneous.

Next, we will verify that $I_\X^\diff$ is saturated.
It suffices to prove the equality 
$I_\X^\diff = I_\X^\diff :_{\Pbar} \mathfrak{M}$.
Let $f\in I_\X^\diff :_{\Pbar} \mathfrak{M}$ and $i\in\{0,\dots,n\}$. 
Then $X_kf \in  I_\X^\diff$ for all $k=0,\dots,n$,
and so $\partial_i(X_kf)\in I_\X$.
Observe that $\partial_i(X_0f)\in I_\X$ and $X_0f \in I_\X$.
In $R$, we have $x_0 \bar{f} =0$. Since $x_0$ is a non-zerodivisor of $R$, 
this yields $\bar{f} =0$, in other words, $f\in  I_\X$.
Consequently, we have 
$$
X_0\partial_i(f) = \partial_i(X_0f) -f\partial_i(X_0) \in I_\X
$$
for all $i=0,\dots,n$.
This also implies $\partial_i f \in I_\X$ for all $i=0,\dots,n$.
Hence $f \in I_\X^\diff$, as was to be shown.

\smallskip\noindent(b)\quad
This follows from Corollary~\ref{Cor-Diff-Powers}.
\end{proof}

When the scheme $\X$ consists of a single monomial point, the homogeneous 
ideal $I_\X^\diff$ can be explicitly described by the following theorem. 
Here the set of all terms in $P$ is given by 
$\mathbb{T}^{n+1}=\{X_0^{\alpha_0}\cdots X_n^{\alpha_n} 
\mid (\alpha_0,\dots,\alpha_n)\in \NN^{n+1} \}$.

\begin{theorem}\label{Prop-Diff-MonomialIdeal}
Let $K$ be a perfect field of characteristic $\charac(K)=p$. 
Assume that the homogeneous vanishing ideal $I_\X$ 
of the scheme $\X$ is a proper monomial ideal of $P$, minimally generated by 
terms $t_1,\dots,t_r \in \mathbb{T}^{n+1}$. 
\begin{enumerate}
	\item[(a)] The ideal $I_\X^\diff$ is a monomial subideal of $I_\X$.
	
	\item[(b)] For $1\le j\le r$, we write 
	$t_j = X_0^{\alpha_{0j}}\cdots X_n^{\alpha_{nj}}$, 
	and for $2\le k\le n+1$ and $1\le j_1\le\cdots\le j_k\le r$, we define
	$t_{j_1,\dots,j_k} := X_0^{\beta_0}\cdots X_n^{\beta_n}$, where, for $i=0,\dots,n$,
	\begin{equation} \label{Equ-DiffElemExponents} 
	\beta_i \;:=\; 
	\begin{cases}
		0 & \text{if $\alpha_{ij_1}=\dots= \alpha_{ij_k}=0$},\\
		\alpha_{ij_1} & 
		\text{if $\alpha_{ij_1}=\dots= \alpha_{ij_k}>0$ and $p\mid \alpha_{ij_1}$},\\
		\alpha_{ij_1}+1 & 
		\text{if $\alpha_{ij_1}=\dots= \alpha_{ij_k}>0$ and $p\nmid \alpha_{ij_1}$},\\
		\max\{\alpha_{ij_1},\dots,\alpha_{ij_k}\} & \text{otherwise}.
	\end{cases}
	\end{equation}
	Then 
	$$
	I_\X^\diff \;=\; \langle\, t_{j_1,\dots,j_k} \mid 
	1\le j_1\le\cdots\le j_k\le r;\; 2\le k\le n+1 \,\rangle.
	$$
\end{enumerate}
\end{theorem}

\begin{proof}
Note that a non-constant term $t = X_0^{\alpha_0}\cdots X_n^{\alpha_n}$
belongs to $I_\X^\diff$ if and only if $\frac{\alpha_it}{X_i} \in I_\X$
for $i=0,\dots,n$. 
	
\smallskip\noindent(a)\quad 
Let $f\in I_\X^\diff$ be a nonzero homogeneous polynomial 
of degree $d$ and write 
$f = c_1t'_1+\cdots + c_mt'_m$ with $c_1,\dots,c_m\in K\setminus\{0\}$ and 
$t'_1,\dots,t'_m\in \mathbb{T}^{n+1}$. It suffices to show that 
$t'_j\in I_\X^\diff$ for $j=1,\dots,m$.
Indeed, let $i\in\{0,\dots,n\}$ and write $t'_j = X_0^{\gamma_{0j}}\cdots X_n^{\gamma_{nj}}$ 
with $(\gamma_{0j},\dots,\gamma_{nj})\in \NN^{n+1}$
and $\gamma_{0j}+\cdots+\gamma_{nj} = d$ for $j=1,\dots,m$.
When $\gamma_{ij}=0$ in $K$, it directly follows that $\partial_i t'_j =0\in I_\X$. 
Now, consider the case, where $\gamma_{ij}\ne 0$ in~$K$. 
In this case, for $1\le j< k\le m$, we have $t'_j \ne t'_k$ and
$$
\partial_i t'_j = \frac{\gamma_{ij}t'_j}{X_i} \ne 
\frac{\gamma_{ik}t'_k}{X_i} = \partial_i t'_k.
$$ 
So, the element $\frac{t'_j}{X_i}$ appears in the 
expansion of $\partial_i f = c_1 \partial_it'_1+\cdots + c_m\partial_it'_m$
with a nonzero coefficient. Since $I_\X$ is a monomial ideal
and $\partial_i f\in I_\X$, this implies $\frac{\partial_i t'_j}{\gamma_{ij}}
=\frac{t'_j}{X_i} \in I_\X$. Consequently, $\partial_i t'_j\in I_\X$
for every $j=1,\dots,m$ and every $i=0,\dots,n$. 
Hence, we find that $t'_1,\dots,t'_m \in I_\X^\diff$.

\smallskip\noindent(b)\quad 
Let 
$$
\Gamma = \{\, t_{j_1,\dots,j_k} \mid 
1\le j_1\le\cdots\le j_k\le r;\; 2\le k\le n+1 \,\}
$$
where $t_{j_1,\dots,j_k} = X_0^{\beta_0}\cdots X_n^{\beta_n}\in \Gamma$ 
has the exponents $\beta_i$ given as in \eqref{Equ-DiffElemExponents}. 
It is easily seen that  $\Gamma \subseteq I_\X$.
When $p=\charac(K)>0$, we also have 
$$
\{ X_l^p\in I_\X \mid l=0,\dots,n\} 
\subseteq \langle \Gamma\rangle\cap I_\X^\diff.
$$ 
Indeed, for $X_l^p\in I_\X$, it is evident that $\partial_i(X_l^p)=0$
for all $i=0,\dots,n$, and hence $X_l^p\in I_\X^\diff$.
Also, there exists $j\in\{1,\dots,r\}$ such that $t_j = X_l^e$
with $1\le e\le p$, and hence $X_l^p$ is a multiple of $t_{j,j}\in\Gamma$. 

Now let $t_{j_1,\dots,j_k} = X_0^{\beta_0}\cdots X_n^{\beta_n}\in \Gamma$. 
We will show that $t_{j_1,\dots,j_k} \in I_\X^\diff$, i.e.,
that $\partial_i t_{j_1,\dots,j_k} \in I_\X$ for all $i=0,\dots,n$.
If $j_1=\cdots=j_k$ then \eqref{Equ-DiffElemExponents} yields
$$
\beta_i = 
\begin{cases}
0 & \text{if $\alpha_{ij_1}=0$},\\
\alpha_{ij_1} &\text{if $\alpha_{ij_1}>0$ and $p\mid \alpha_{ij_1}$,}\\
\alpha_{ij_1}+1 &\text{if $\alpha_{ij_1}>0$ and $p\nmid \alpha_{ij_1}$}.
\end{cases}
$$
More specifically, if $p\mid \alpha_{ij_1}$, then $\beta_i=0$ in $K$,  
and $\beta_i =\alpha_{ij_1}+1$ otherwise. This directly implies that 
$\partial_i t_{j_1,\dots,j_k}=\beta_i\frac{t_{j_1,\dots,j_l}}{X_i}\in I_\X$ 
for all $i=0,\dots,n$,
and consequently $t_{j_1,\dots,j_k}\in I_\X^\diff$.

Next, let us assume that $j_{k_1}\ne j_{k_2}$ for some $k_1,k_2 \in\{1,\dots,k\}$
and let $i\in\{0,\dots,n\}$.
Obviously, $p\mid \beta_i$ implies $\partial_i t_{j_1,\dots,j_k}=0\in I_\X$.
So, we assume that $p\nmid \beta_i$. By~\eqref{Equ-DiffElemExponents}, we see that
$\beta_i \ge \alpha_{ij_l}$ for all $l\in\{1,\dots,k\}$ and 
for all $i\in \{0,\dots,n\}$. Consider the following two cases:
\begin{itemize}
	\item $\alpha_{ij_1}=\dots= \alpha_{ij_k}>0$:\
	Then $\beta_i = \alpha_{ij_1}+1$ and 
	$$
	\partial_i (t_{j_1,\dots,j_k}) =   \frac{\beta_it_{j_1,\dots,j_k}}{X_i} 
	= \beta_iX_0^{\beta_0}\cdots X_{i-1}^{\beta_{i-1}} X_i^{\beta_i-1}
	X_{i+1}^{\beta_{i+1}}\cdots X_n^{\beta_n}
	$$ 
	is a multiple of $t_{j_1}$, 
	and hence it belongs to  $I_\X$.
	
	\item There exist $l_1,l_2\in \{1,\dots,k\}$ such that 
	$\alpha_{ij_{l_1}}>\alpha_{ij_{l_2}}$: \
	Without loss of generality, we may assume that $l_2=1$. By~\eqref{Equ-DiffElemExponents}, we have  
	$$
	\beta_i = \max\{\alpha_{ij_1},\dots,\alpha_{ij_k}\}
	\ge \alpha_{ij_{l_1}}\ge \alpha_{ij_1}+1.
	$$ 
	Thus, $\partial_i (t_{j_1,\dots,j_k})$ is 
	a multiple of $t_{j_1}$, and so it is an element of $I_\X$.
\end{itemize}
Hence, we get $\partial_i (t_{j_1,\dots,j_l}) \in I_\X$ 
for all $i\in\{0,\dots,n\}$, which means 
$t_{j_1,\dots,j_l}\in I_\X^\diff$.
Consequently, $\langle \Gamma \rangle\subseteq I_\X^\diff$.

For the ``$\supseteq$'' inclusion, 
let $u = X_0^{\gamma_0} \cdots X_n^{\gamma_n}$ with $(\gamma_0, \dots, \gamma_n)
\in \NN^{n+1}$  be a term in $I_\X^\diff$. By \eqref{Equ-DiffIdeal}, 
we have $u\subseteq I_\X$ and $\partial_iu\in  I_\X$ for $i=0,\dots,n$.
Since $\{t_1,\dots,t_r\}$ is the minimal monomial system of generators of $I_\X$,
there exists $j\in\{1,\dots,r\}$ such that $t_j \mid u$.
Let $T_u$ be the set of all terms in $\{t_1,\dots,t_r\}$ which divide $u$.
After renumbering the indices of $t_j$, we may assume that 
$T_u =\{ t_1,\dots, t_s\}$ with $s\le r$. 
Obviously, we have 
\begin{equation}\label{Equ-DiffExponents}
	\gamma_i \ge \max\{\alpha_{ij} \mid j=1,\dots,s\} \qquad (i=0,\dots,n).
\end{equation}
Now we distinguish between two cases 
concerning the characteristic $p$ of $K$. 

\smallskip\noindent \textbf{Case 1:} $p=0$. 
If $s=1$, then $t_1 \mid u$ and $t_j\nmid u$ for $j=2,\dots,r$,
in particular, $t_j \nmid \partial_i u$ for all $i=0,\dots,n$ and $j=2,\dots,r$.
If $\alpha_{i1}>0$, then $\partial_i u \in I_\X$ implies $t_1\mid \partial_i u$,
and hence $\gamma_i\ge \alpha_{i1}+1$ for $i=0,\dots,n$.
It follows that $t_{1,1} \mid u$ or $u\in \langle \Gamma \rangle$.

Suppose that $s\ge 2$ and consider the following cases:
\begin{itemize}
	\item $s\le n+1$: \ Let $j_1=1,\dots,j_s=s$ and 
	$t_{j_1,\dots,j_s} = X_0^{\beta_0} \cdots X_n^{\beta_n}$,
	where the $\beta_i$ are given as in \eqref{Equ-DiffElemExponents}, 
	i.e., letting $k=s$ we have
	\begin{equation} \label{Equ-DiffElemExponents 2} 
		\beta_i \;=\; 
		\begin{cases}
			0 & \text{if $\alpha_{ij_1}=\dots= \alpha_{ij_k}=0$},\\
			\alpha_{ij_1}+1 & 
			\text{if $\alpha_{ij_1}=\dots= \alpha_{ij_k}>0$},\\
			\max\{\alpha_{ij_1},\dots,\alpha_{ij_k}\} & \text{otherwise}.
		\end{cases}
	\end{equation}
	We will prove that $t_{j_1,\dots,j_s}$ divides $u$.
	Consider any $i\in\{0,\dots,n\}$.
	\begin{enumerate}
		\item[$\ast$] 
		If $\alpha_{i1}=\dots=\alpha_{is}=0$ then $\gamma_i\ge \beta_i=0$.
		
		\item[$\ast$] 
		If $\alpha_{i1}=\dots=\alpha_{is}>0$, then $\gamma_i\ge \alpha_{i1}$.
		When $\gamma_i = \alpha_{i1}$, we would have $t_j \nmid \partial_i u$ 
		for all $j=1,\dots,s$, and hence $\partial_i u\notin I_\X$, which is impossible.
		Thus, we must have $\gamma_i\ge \alpha_{i1}+1=\beta_i$.
		
		\item[$\ast$] 
		Suppose there exist $j,k\in\{1,\dots,s\}$ such that 
		$\alpha_{ij}<\alpha_{ik}$. In this case, by \eqref{Equ-DiffExponents}, 
		we also have $\gamma_i\ge \max\{\alpha_{ij} \mid j=1,\dots,s\} = \beta_i$.
	\end{enumerate}
	Hence, $\gamma_i \ge \beta_i$ for all $i\in\{0,\dots,n\}$,
	and subsequently $u$ is divisible by $t_{j_1,\dots,j_s}$.
	Therefore, we get  $u\in \langle \Gamma \rangle$.
	
	\item $s> n+1$:\quad Let $\{i_{k+1},\dots,i_{n+1}\}$ be the set of indices
	$i \in \{0,\dots,n\}$ for which $\alpha_{i1}=\dots=\alpha_{is}$.
	The complement of this set within $\{0,\dots,n\}$ is written as
	$\{i_1,\dots,i_k\}$.
	For each $l\in \{1,\dots,k\}$, we recursively choose an index 
	$j_l \in \{1,\dots,s\}$ such that
	$$
	\alpha_{i_l j_l} = \min\{\, \alpha_{i_l j} \mid 
	j\in \{1,\dots,s\}\setminus\{j_1,\dots,j_{l-1}\} \,\}.
	$$
	Let $t_{j_1,\dots,j_k}= X_0^{\beta_0} \cdots X_n^{\beta_n}$, where 
	the $\beta_i$ are given as in \eqref{Equ-DiffElemExponents 2}.
	We will demonstrate that $t_{j_1,\dots,j_k}$ divides $u$.
	Note that $\{0,\dots,n\} = \{i_1,\dots,i_{n+1}\}$.
	\begin{itemize}
		\item[$\ast$] 
		For $l\in \{k+1,\dots,n+1\}$, if $\alpha_{i_l 1}=\dots=\alpha_{i_l s}=0$
		then $\gamma_{i_l}\ge 0=\beta_{i_l}$; and
		if $\alpha_{i_l1}=\dots=\alpha_{i_ls}>0$ then $\gamma_{i_l}\ge \alpha_{i_l1}+1=\beta_{i_l}$, 
		since otherwise $\gamma_{i_l}=\alpha_{i_l\,1}$ 
		would imply that $t_e \nmid \partial_{i_l} u$ 
		for all $e=1,\dots,s$, and so $\partial_{i_l}u\notin I_\X$, a contradiction.
		
		\item[$\ast$] 
		For $l\in\{1,\dots,k\}$,
		if $\alpha_{i_lj_1}=\dots=\alpha_{i_lj_k}$ then 
		$$
		\alpha_{i_lj_1} = \min\{\, \alpha_{i_lj} \mid 
		j\in \{1,\dots,s\}\}
		\quad\text{and}\quad  
		\beta_{i_l} =\alpha_{i_lj_1}+1.
		$$
		Due to our choice of $i_l$ and \eqref{Equ-DiffExponents}, 
		we have 
		$$
		\gamma_{i_l} \ge \max\{\alpha_{i_lj} \mid j=1,\dots,s\}
		\ge \alpha_{i_lj_1}+ 1 =\beta_{i_l}.
		$$
		If $\alpha_{i_lj_l}\ne \alpha_{i_lj_e}$ for some $e\in\{1,\dots,k\}$, then 
		we also have
		$$
		\gamma_{i_l} \ge \max\{\alpha_{i_lj} \mid j=1,\dots,s\}
		\ge \max\{\alpha_{i_lj_e} \mid e=1,\dots,k\} = \beta_{i_l}.
		$$
	\end{itemize}
	So, $\gamma_{i_l} \ge \beta_{i_l}$ for all $l=1,\dots,n+1$,
	and hence $t_{j_1,\dots,j_k} \mid u$.
	Thus, $u\in \langle \Gamma \rangle$.
\end{itemize}
Consequently, for $p=0$, we have established the inclusion
$\langle \Gamma \rangle \supseteq I_\X^\diff$.

\smallskip\noindent \textbf{Case 2:} $p>0$. 
Our aim here is to show that the term
$u = X_0^{\gamma_0} \cdots X_n^{\gamma_n}$ in~$I_\X^\diff$ 
also belongs to~$\langle \Gamma\rangle$.
Recall that $T_u=\{t_1,\dots,t_s\}$ is defined as the set of 
all terms in $\{t_1,\dots,t_s\}$ that divide $u$.
We proceed by examining the following subcases:
\begin{itemize}
	\item $s=1$:\quad We have $t_1 \mid u$ and $t_j\nmid u$ for $j=2,\dots,r$.
	Obviously, $\gamma_i\ge \alpha_{i1}$ for all $i=0,\dots,n$.
	If $p\nmid \gamma_i$, then $t_j \nmid \partial_i u$ for all $j=2,\dots,r$.
	From $\partial_i u \in I_\X\setminus\{0\}$ it follows that 
	$t_1$ divides $\partial_i u$, and so $\gamma_i\ge \alpha_{i1}+1$.
	If $p\mid \gamma_i$ and $p\nmid \alpha_{i1}$, then $\gamma_i\ge \alpha_{i1}+1$.
	By \eqref{Equ-DiffElemExponents}, we find $t_{1,1} \mid u$ 
	or $u\in \langle \Gamma \rangle$.
	
	\item $2\le s\le n+1$: \ Let $j_1=1,\dots,j_s=s$ and 
	$t_{j_1,\dots,j_s} = X_0^{\beta_0} \cdots X_n^{\beta_n}$,
	where the $\beta_i$ are given as in \eqref{Equ-DiffElemExponents}.
	We will verify that $t_{j_1,\dots,j_s}$ divides $u$.
	Let $i\in\{0,\dots,n\}$.
	\begin{enumerate}
		\item[$\ast$] 
		If $\alpha_{i1}=\dots=\alpha_{is}=0$ then $\gamma_i\ge \beta_i=0$.
		
		\item[$\ast$] 
		If $\alpha_{i1}=\dots=\alpha_{is}>0$ and $p\mid \alpha_{i1}$, then 
		$\gamma_i\ge \alpha_{i1}=\beta_i$.
		
		\item[$\ast$] 
		If $\alpha_{i1}=\dots=\alpha_{is}>0$ and $p\nmid \alpha_{i1}$, 
		then $\gamma_i\ge \alpha_{i1}$.
		When $\gamma_i = \alpha_{i1}$, we would have $t_j \nmid \partial_i u$ 
		for all $j=1,\dots,s$, and hence $\partial_i u\notin I_\X$, which is impossible.
		Thus, we must have $\gamma_i\ge \alpha_{i1}+1=\beta_i$.
		
		\item[$\ast$] 
		Suppose there exist $j,k\in\{1,\dots,s\}$ such that 
		$\alpha_{ij}<\alpha_{ik}$. In this case, by \eqref{Equ-DiffExponents}, 
		we also have $\gamma_i\ge \max\{\alpha_{ij} \mid j=1,\dots,s\} = \beta_i$.
	\end{enumerate}
	Accordingly, $\gamma_i \ge \beta_i$ for all $i\in\{0,\dots,n\}$,
	which means $u$ is divisible by $t_{j_1,\dots,j_s}$.
	Hence  $u\in \langle \Gamma \rangle$.

	\item $s> n+1$:\quad Let $\{i_{k+1},\dots,i_{n+1}\}$ be the set of indices
	$i \in \{0,\dots,n\}$ for which $\alpha_{i1}=\dots=\alpha_{is}$ and 
	$$
	\{i_1,\dots,i_k\}=\{0,\dots,n\}\setminus\{i_{k+1},\dots,i_{n+1}\}.
	$$
	For each $l\in \{1,\dots,k\}$, we recursively choose an index 
	$j_l \in \{1,\dots,s\}$ such that
	$$
	\alpha_{i_l j_l} = \min\{\, \alpha_{i_l j} \mid 
	j\in \{1,\dots,s\}\setminus\{j_1,\dots,j_{l-1}\} \,\}.
	$$
	Let $t_{j_1,\dots,j_k}= X_0^{\beta_0} \cdots X_n^{\beta_n}$, where 
	the $\beta_i$ are given as in \eqref{Equ-DiffElemExponents}.
	We will check that $t_{j_1,\dots,j_k}\mid u$.
	\begin{itemize}
		\item[$\ast$] 
		If $l\in \{k+1,\dots,n+1\}$ and $\alpha_{i_l\,1}=\dots=\alpha_{i_l\,s}=0$,
		then $\gamma_{i_l}\ge 0=\beta_{i_l}$.
		
		\item[$\ast$] 
		If $l\in \{k+1,\dots,n+1\}$ and $\alpha_{i_l\,1}=\dots=\alpha_{i_l\,s}>0$
		and $p\mid \alpha_{i_l\,1}$, then $\gamma_{i_l}\ge \alpha_{i_l\,1}=\beta_{i_l}$.
		
		\item[$\ast$] 
		If $l\in \{k+1,\dots,n+1\}$ and $\alpha_{i_l\,1}=\dots=\alpha_{i_l\,s}>0$
		and $p\nmid \alpha_{i_l\,1}$, then $\gamma_{i_l}\ge \alpha_{i_l1}+1=\beta_{i_l}$, 
		since otherwise $\gamma_{i_l}=\alpha_{i_l\,1}$ 
		would imply that $t_e \nmid \partial_{i_l} u$ 
		for all $e=1,\dots,s$, and so $\partial_{i_l}u\notin I_\X$, a contradiction.
		
		\item[$\ast$] 
		If $l\in\{1,\dots,k\}$ and $\alpha_{i_lj_1}=\dots=\alpha_{i_lj_k}$,	then 
		$$
		\alpha_{i_lj_1} = \min\{\, \alpha_{i_lj} \mid 
		j\in \{1,\dots,s\}\}
		\quad\text{and}\quad  
		\beta_{i_l} =\alpha_{i_lj_1}+1.
		$$
		By the choice of $i_l$ and \eqref{Equ-DiffExponents}, 
		we have 
		$$
		\gamma_{i_l} \ge \max\{\alpha_{i_lj} \mid j=1,\dots,s\}
		\ge \alpha_{i_lj_1}+ 1 =\beta_{i_l}.
		$$
		
		\item[$\ast$] 
		If $l\in\{1,\dots,k\}$ and $\alpha_{i_lj_l}\ne \alpha_{i_lj_e}$ for some $e\in\{1,\dots,k\}$, then 
		$$
		\gamma_{i_l} \ge \max\{\alpha_{i_lj} \mid j=1,\dots,s\}
		\ge \max\{\alpha_{i_lj_e} \mid e=1,\dots,k\} = \beta_{i_l}.
		$$
	\end{itemize}
	Thus, $\gamma_{i_l} \ge \beta_{i_l}$ for all $l=1,\dots,n+1$,
	and so $t_{j_1,\dots,j_k} \mid u$.
	Hence $u\in \langle \Gamma \rangle$.
\end{itemize}
Altogether, the inclusion $\langle \Gamma\rangle \supseteq I_\X^\diff$
also holds true in the case of positive characteristic, 
thereby completing the proof.
\end{proof}

\begin{remark} 
The proof of Proposition~\ref{Prop-Diff-MonomialIdeal} remains valid for 
any arbitrary nonzero monomial ideal $I$ in $P$. 
In particular, a finite monomial systems of generators for $I^\diff$
can be explicitly constructed from the minimal generators of $I$.
\end{remark}

\begin{corollary}\label{Cor-Diff-MonomialIdeal}
Let $K$ be a perfect field of characteristic $\charac(K)=p$. 	
Suppose that $\X$ contains exactly one point $p=(1:0:\dots:0)$. 
	\begin{enumerate}
		\item[(a)] If $I_\X = \langle X_1,\dots,X_n\rangle^m$ for some $m\ge 1$, 
		then 
		$$
		I_\X^\diff = 
		\begin{cases}
			\langle X_1,\dots,X_n\rangle^{m+1}& \text{ if $p\nmid m$,}\\
			\langle X_1,\dots,X_n\rangle^{m+1}+
			\big\langle X_1^{p\alpha_1}\cdots X_n^{p\alpha_n} \mid 
			\sum_{i=1}^n\alpha_i=\frac{m}{p} \big\rangle
			& \text{ if $p\mid m$.}
		\end{cases}
		$$
		
		\item[(b)] If $I_\X = \langle X_1^{k_1},\dots,X_n^{k_n}\rangle$ for 
		some integers $k_1,\dots,k_n\ge 1$, then 
		$$
		I_\X^\diff = 
		\begin{cases}
			\langle X_1^{k_1+1},\dots,X_n^{k_n+1}\rangle
			+ \langle X_i^{k_i}X_j^{k_j} \mid 1\le i<j\le n\rangle 
			& \text{ if $p=0$,}\\
			\begin{array}{l}
			\langle X_i^{k_i}\mid i=1,\dots,n;\, p \mid k_i\rangle +
			\langle X_i^{k_i+1}\mid i=1,\dots,n;\, p\nmid k_i\rangle \\
			+ \langle X_i^{k_i}X_j^{k_j} \mid 1\le i<j\le n\rangle
			\end{array}
			& \text{ if $p>0$.}
		\end{cases}
		$$
	\end{enumerate}
\end{corollary}

\begin{proof}
Let $\{t_1,\dots,t_r\}$ be the minimal monomial system of generators of $I_\X$. 
Define $\Gamma$ as in the proof of Theorem~\ref{Prop-Diff-MonomialIdeal}(b):
$$
\Gamma = \{ t_{j_1,\dots,j_k} \mid 1\le j_1\le\cdots\le j_k\le r;\, 
2\le k\le n+1\}.
$$ 
Then $I_\X^\diff = \langle \Gamma \rangle$. 
In both considered cases, no term in $\Gamma$ is divisible by $X_0$.

\smallskip\noindent(a)\quad We have $r=\binom{n+m-1}{m}$ and 
$t_j =X_1^{\alpha_{1j}}\cdots X_n^{\alpha_{nj}}$
with $\alpha_{1j}+\cdots+\alpha_{nj}=m$ for $j=1,\dots,r$.
For $1\le j_1< j_2\le r$, there always exists $i\in\{1,\dots,n\}$ such that
$\alpha_{ij_1} \ne \alpha_{ij_2}$, and so $\deg(t_{j_1,j_2})\ge m+1$. 
Also, for $i\in\{1,\dots,n\}$, let $u_{ij} = X_it_j$ and 
$\{i_1,\dots,i_k\} = \{i \mid \alpha_{ij}\ne 0\}$, and write
$t_{j_l} = \frac{u_{ij}}{X_{i_l}}$ for $l=1,\dots,k$. 
By \eqref{Equ-DiffElemExponents}, we find $u_{ij} = t_{j_1,\dots,j_{k}}\in \Gamma$.
Thus 
$$
\langle X_1,\dots,X_n\rangle^{m+1} \subseteq \langle \Gamma \rangle = I_\X^\diff.
$$
Let us consider the following cases:
\begin{itemize}
	\item $p\nmid m$:\quad 
	Clearly, $X_l^m\notin \Gamma$ for $l=0,\dots,n$. 
	If $t_j =X_1^{\alpha_{1j}}\cdots X_n^{\alpha_{nj}}\in\{X_l^m \mid l=0,\dots,n\}$, then $p \mid m=\sum_{i=1}^n \alpha_{ij}$ implies $p\nmid \alpha_{ij}$
	for some $i\in\{1,\dots,n\}$. For such $i$, it follows that 
	$\partial_i (t_j) = \alpha_{ij}\frac{t_j}{X_i} \notin I_\X$, 
	and hence $t_j\notin I_\X^\diff$, in particular, $t_j\notin \Gamma$.
	Hence, no term of degree $\le m$ belongs to $\Gamma$.
	This implies $\langle X_1,\dots,X_n\rangle^{m+1} =\langle \Gamma \rangle$.
	
	\item $p\mid m$:\quad 
	For every $(\alpha_1,\dots,\alpha_n)\in\NN^n$ with $\sum_{i=1}^n p\alpha_i =m$, 
	we have $X_1^{p\alpha_1}\cdots X_n^{p\alpha_n}\in I_\X^\diff$, and so
	$X_1^{p\alpha_1}\cdots X_n^{p\alpha_n} \in \Gamma$.
	If $t_j =X_1^{\alpha_{1j}}\cdots X_n^{\alpha_{nj}}$ is not in 
	$\{ X_1^{p\alpha_1}\cdots X_n^{p\alpha_n} \mid 
	\sum_{i=1}^n\alpha_i=\frac{m}{p}\}$, then we find 
	$p\nmid \alpha_{ij}$ for some $i\in\{1,\dots,n\}$.
	Thus, $\partial_i (t_j) = \alpha_{ij}\frac{t_j}{X_i} \notin I_\X$, 
	and hence $t_j\notin I_\X^\diff$. This implies that $t_j\notin \Gamma$. 
	Therefore, we get the equality
	$$
	I_\X^\diff = 
		\langle X_1,\dots,X_n\rangle^{m+1}+
		\big\langle X_1^{p\alpha_1}\cdots X_n^{p\alpha_n} \mid 
		\textstyle{\sum_{i=1}^n\alpha_i=\frac{m}{p}} \big\rangle.
	$$
\end{itemize}

\smallskip\noindent(b)\quad 
We have $t_1= X_1^{k_1},\dots,t_n=X_n^{k_n}$ and $r=n$.
Since $t_1,\dots,t_n$ are pairwise coprime, by \eqref{Equ-DiffElemExponents},
any other term in $\Gamma$ is a multiple of $t_{i,j}$ for some $1\le i\le j\le n$.
Hence $I_\X^\diff = \langle \Gamma \rangle =
\langle t_{i,j} \mid 1\le i\le j\le n\rangle$.
It is enough to describe the set $\{t_{i,j} \mid 1\le i\le j\le n\}$.
For this end, we distinguish between the following cases.
\begin{itemize}
	\item $p=0$:\quad 
	It is easily seen that $t_{i,i} = X_i^{k_i+1}$ for
	$i\in\{1,\dots,n\}$ and $t_{i,j} = X_i^{k_i}X_j^{k_j}$ 
	for all $1\le i<j\le n$. Thus, 
	$$
	\{ t_{i,j} \mid 1\le i\le j\le n \} \;=\;
	  \{ X_1^{k_1+1},\dots,X_n^{k_n+1} \}
	  \cup \{ X_i^{k_i}X_j^{k_j} \mid 1\le i<j\le n \}.
	$$
	
	\item $p>0$:\quad 
	Observe that if $p\nmid k_i$ then $t_{i,i} = X_i^{k_i+1} \in \Gamma$;
	and if $p\mid k_i$ then $t_{i,i} = X_i^{k_i} \in \Gamma$.
	Also, $t_{i,j} = X_i^{k_i}X_j^{k_j}\in \Gamma$ for all $1\le i<j\le n$.
	Hence 
	\begin{align*}
	\{ t_{i,j} \mid 1\le i\le j\le n \} \; =&\;
	\{X_i^{k_i} \mid i=1,\dots,n;\ p\mid k_i \} 
	\cup  \{X_i^{k_i+1} \mid i=1,\dots,n;\ p\nmid k_i \} \\
	&\cup  \{ X_i^{k_i}X_j^{k_j} \mid 1\le i<j\le n \}.
	\end{align*}
\end{itemize}
Therefore, claim (b) is fully demonstrated.
\end{proof}

\begin{remark}\label{Rem-TheFirstFattening}
For $K$-rational points $p_1,\dots,p_s$ and positive integers $m_1,\dots, m_s$, 
the saturated homogeneous ideal $I_\X := I_{p_1}^{m_1}\cap\cdots\cap I_{p_s}^{m_s}$
of $P$ defines a fat point scheme $\X= m_1p_1+\cdots+m_sp_s$ in $\mathbb{P}^n_K$.
The first fattening of this scheme is the fat point scheme 
$\Y=(m_1+1)p_1+\cdots+(m_s+1)p_s$ in~$\mathbb{P}^n_K$. 
When $\X= m_1p_1+\cdots+m_sp_s$ and $\charac(K)\nmid m_j$ for every $j=1,\dots,s$, 
Proposition~\ref{Prop-Diff-HomogVanishingIdeal} 
and Corollary~\ref{Cor-Diff-MonomialIdeal}(a) imply that $I_\X^\diff$ is
the homogeneous vanishing ideal of the first fattening $\Y$ of $\X$ in $\mathbb{P}^n_K$.	
\end{remark}

\begin{definition}
Let $\X$ be a 0-dimensional scheme in $\mathbb{P}^n_K$ with homogeneous 
vanishing ideal~$I_\X$. Then the saturated homogeneous ideal $I_\X^\diff$
defines a 0-dimensional scheme $\Y$ in $\mathbb{P}^n_K$.
The scheme $\Y$ is called the \textbf{differential fattening} of $\X$.
\end{definition}

With this terminology in hand, we can now provide an affirmative answer to 
Question~\ref{Question 1}, regarding the existence of 
the canonical exact sequence of the K\"ahler differential module $\Omega^1_R$ 
for an arbitrary 0-dimensional scheme $\X$.

\begin{theorem}\label{Thm-CanonicalExactSequence}
Let $\X$ be a 0-dimensional scheme in $\mathbb{P}^n_K$ with homogeneous 
vanishing ideal~$I_\X$, and let $\Y$ be the differential fattening of $\X$
with its homogeneous vanishing ideal~$I_\Y$.
The sequence of graded $R$-modules
$$
0 \longrightarrow I_\X/ I_\Y \stackrel{\delta}{\longrightarrow}
R^{n+1}(-1)  \stackrel{\gamma}{\longrightarrow} \Omega^1_R
\longrightarrow 0
$$
is exact, where $\delta(f + I_\Y) 
=(\partial_0(f),\dots, \partial_n(f)) 
= \sum_{i=0}^{n}\partial_i(f)e_i$ 
for every homogeneous polynomial $f\in I_\X$ and
$\gamma(e_i) = dx_i$ for $i=0,\dots,n$.
\end{theorem}

\begin{proof}
Notice that $R^{n+1}(-1)\cong \Omega^1_{P}/I_\X\Omega^1_{P}$
and $I_\Y= I_\X^\diff$.
It suffices to prove that the sequence of graded $R$-module
$$
0 \longrightarrow I_\X/ I_\X^\diff \stackrel{\delta^*}{\longrightarrow}
\Omega^1_{P}/I_\X\Omega^1_{P}  \stackrel{\gamma^*}{\longrightarrow} 
\Omega^1_R \longrightarrow 0
$$
is exact, where $\delta^*(f + I_\X^\diff) =df + I_\X\Omega^1_{P}$ 
for each homogeneous polynomial $f\in I_\X$ and
$\gamma^*(dX_i + I_\X\Omega^1_{P}) = dx_i$ for $i=0,\dots,n$.
For this, we will first show that $\delta^*$ is well-defined.
In fact, let $f\in I_\X\setminus\{0\}$ be such that 
$f+  I_\X^\diff = 0+I_\X^\diff$,
i.e., that $f\in I_\X^\diff$.
By definition of $I_\X^\diff$, we have $\partial_i(f) \in I_\X\Omega^1_{P}$
for all $i=0,\dots,n$, and hence 
$$
\delta^*(f+I_\X^\diff) = df +I_\X\Omega^1_{P} 
=\sum_{i=0}^n \partial_i(f)dX_i +I_\X\Omega^1_{P} = 0 +I_\X\Omega^1_{P}.
$$ 
This means that $\delta^*$ is well-defined.
Also, it is $R$-linear map. In fact, for $f_1,f_2\in I_\X$ and $g_1,g_2\in P$,
we have
\begin{align*}
	\delta^*(f_1g_2+f_2g_2 + I_\X^\diff) 
	&= d(f_1g_2+f_2g_2) + I_\X\Omega^1_{P}
	= (g_1df_1 +g_2df_2)+ I_\X\Omega^1_{P}\\
	&= (g_1+I_\X)\cdot \delta^*(f_1 + I_\X^\diff)
	+ (g_2+I_\X)\cdot \delta^*(f_2 + I_\X^\diff).
\end{align*}
	
Clearly, both $\delta^*$ and $\gamma^*$ are homogeneous of degree 0, 
and $\gamma^*$ is surjective. Moreover, by \cite[Proposition 4.12]{Kun1986}, 
we have $\Omega^1_R \cong \Omega^1_{P}/(dI_\X + I_\X \Omega^1_{P})$
and $\Im(\delta^*) = (dI_\X + I_\X \Omega^1_{P})/I_\X \Omega^1_{P}$,
and hence we get $\Im(\delta^*) =\Ker(\gamma^*)$.
	
Next, it remains to show that $\delta^*$ is injective. Indeed, let $f\in I_\X$
be such that $\delta^*(f+ I_\X^\diff)=0$. Then 
$df + I_\X \Omega^1_{P} = 0 + I_\X \Omega^1_{P}$, and hence
$df = \sum_{i=0}^n \partial_i(f)dX_i\in I_\X \Omega^1_{P}$. 
It follows that $\partial_i(f) \in I_\X$ for all $i=0,\dots,n$.
Therefore $f\in I_\X^\diff$, and hence 
$f+ I_\X^\diff =0+I_\X^\diff$, as desired.
\end{proof}

\begin{remark}\label{Rem-CaExaxtSeq-FatPointSch}
Theorem~\ref{Thm-CanonicalExactSequence} and Remark~\ref{Rem-TheFirstFattening}
together generalize the result of \cite[Theorem 1.7]{KLL2015} 
for fat point schemes $\X=m_1p_1+\cdots+m_sp_s$ in $\mathbb{P}^n_K$.
This holds even when the base field $K$ is not necessary of characteristic zero,
provided that $\charac(K)\nmid m_j$ for all $j=1,\dots,s$.
\end{remark}

\begin{corollary}\label{Cor-HFOmega^1_R}
Let $\X$ be a 0-dimensional scheme in $\mathbb{P}^n_K$, and let
$\Y$ be the differential fattening of $\X$. 
\begin{enumerate}
	\item[(a)] For $i\in \NN$, we have 
	$$
	\HF_{\Omega^1_R}(i) = (n+1)\HF_\X(i-1) +\HF_\X(i) - \HF_\Y(i).
	$$
	In particular, $\HP(\Omega^1_R)=(n+2)\deg(\X) -\deg(\Y)$
	and the regularity index of $\Omega^1_R$ is bounded by
	$\operatorname{ri}(\Omega^1_R)\le \max\{r_\X+1, r_\Y\}$.
	
	\item[(b)] 
	Suppose $\X=m_1p_1+\cdots+m_sp_s$ is a fat point scheme in $\mathbb{P}^n_K$ 
	and let $\charac(K)=p$.
	
	\begin{enumerate}
		\item[(i)] We have 
		\begin{align*}
			\HP(\Omega^1_R) \;=\; 
			& \sum_{j=1}^s \Big[ (n+2)\binom{m_j+n-1}{n} -\binom{m_j+n}{n} \Big] \\
			& +
			\sum_{p\mid m_j}\#\Big\{(\alpha_1,\dots,\alpha_n)\in\NN^n \mid 
			\sum_{i=1}^{n}\alpha_i = \frac{m_j}{p}  \Big\}.
		\end{align*}
		
		\item[(ii)] We have $m_1=\cdots=m_s=1$ if and only if $\HP(\Omega^1_R)=s$.
	\end{enumerate}
\end{enumerate}
\end{corollary}

\begin{proof}
Claim (a) directly follows from Theorem~\ref{Thm-CanonicalExactSequence}.
Now we will prove two parts in claim (b). For (i), we have 
$\HP(\Omega^1_R) \;=\; (n+2)\deg(\X) - \deg(\Y)$ by (a).
For $j\in\{1,\dots,s\}$, let $\X_j = m_jp_j$ and let 
$\Y_j$ be the differential fattening of $\X_j$.
Then $\deg(\X)=\sum_{j=1}^s \deg(\X_j) =  \sum_{j=1}^s\binom{m_j+n-1}{n}$
and $\deg(\Y) = \sum_{j=1}^s \deg(\Y_j)$.
So, it is sufficient to locally examine $\deg(\Y_j)$. 
After a homogeneous linear change of coordinates, we may assume $p_j=(1:0:\cdots:0)$. 
Then $I_{\X_j} = \langle X_1,\dots,X_n\rangle^{m_j}$.
If $p\nmid m_j$, then, by Corollary~\ref{Cor-Diff-MonomialIdeal}(a),
$I_{\Y_j} = \langle X_1,\dots,X_n\rangle^{m_j+1}$ and $\Y_j=(m_j+1)p_j$, 
and consequently $\deg(\Y_j) = \binom{m_j+n}{n}$.
When $p\mid m_j$, Corollary~\ref{Cor-Diff-MonomialIdeal}(a) yields 
$$
I_{\Y_j} = \langle X_1,\dots,X_n\rangle^{m_j+1}+
\big\langle X_1^{p\alpha_1}\cdots X_n^{p\alpha_n} \mid 
\textstyle{\sum_{i=1}^n\alpha_i=\frac{m_j}{p}} \big\rangle.
$$
Letting $A=K[X_1,\dots,X_n]$, it follows that $I_{\Y_j}^\deh = I_{\Y_j}\cap A$ and
$$
\deg(\Y_j) = \dim_K(A/(I_{\Y_j}\cap A))
= \binom{m_j+n}{n} -\#\Big\{(\alpha_1,\dots,\alpha_n)\in\NN^n \mid 
\sum_{i=1}^{n}\alpha_i = \frac{m_j}{p}  \Big\}.
$$
Hence we get the desired formula for $\HP(\Omega^1_R)$ as in (i). 

Finally, from (i), we have 
\begin{align*}
\HP(\Omega^1_R)=
&\sum_{j=1}^s \Big[ (n+2)\binom{m_j+n-1}{n} -\binom{m_j+n}{n} \Big] \\
& +\sum_{p\mid m_j}\#\Big\{(\alpha_1,\dots,\alpha_n)\in\NN^n \mid 
\sum_{i=1}^{n}\alpha_i = \frac{m_j}{p}  \Big\}\\
=& \sum_{j=1}^s \Big[ \binom{m_j+n-1}{n} +\binom{m_j+n-1}{n-1}(m_j-1) \Big]\\
& +\sum_{p\mid m_j}\#\Big\{(\alpha_1,\dots,\alpha_n)\in\NN^n \mid 
\sum_{i=1}^{n}\alpha_i = \frac{m_j}{p}  \Big\}.
\end{align*}
Therefore, $\HP(\Omega^1_R)=s$ if and only if $m_1=\cdots=m_s=1$,
and (ii) follows.
\end{proof}

In the sense of \cite[Definition 4.3.1]{KR2016}, a 0-dimensional scheme $\X$
is termed a curvilinear scheme if, after a homogeneous linear change of coordinates,
the ideal $I_{p_j}$ for each $j\in\{1,\dots,s\}$ can be expressed as 
$I_{p_j} = \langle X_1,\dots,X_{n-1}, X_n^{k_j} \rangle$
for some positive integer~$k_j$. 
For such a curvilinear scheme $\X$, the ideal $J_{p_j} = I_{p_j}^\deh$
in $A=K[X_1,\dots,X_n]$ is also generated by $\{X_1,\dots,X_{n-1}, X_n^{k_j}\}$ 
after a linear change of coordinates. 
The affine ideal of $\X$ in $D_+(X_0)\cong\mathbb{A}^n_K$ is 
$J_\X = J_{p_1}\cap\cdots\cap J_{p_s}$.
By the Chinese Remainder Theorem, the affine coordinate ring $S=A/J_\X$ 
satisfies $S = \mathcal{O}_1\times\cdots\times\mathcal{O}_s$, where 
$\mathcal{O}_j\cong A/J_{p_j}$ is the 0-dimensional local ring of~$\X$ at~$p_j$.
The next corollary describes the Hilbert polynomial of the K\"ahler differential 
module $\Omega^1_R$ for a curvilinear scheme 
(see also \cite[Corollary 6.7]{KLL2025} for an alternative proof).

\begin{corollary} \label{Cor-Curvilinear}
Let $\X$ be a 0-dimensional curvilinear scheme in $\mathbb{P}^n_K$,
and let $k_j =\dim_K(\mathcal{O}_j)=\dim_K(A/J_{p_j})$ for each $j\in\{1,\dots,s\}$.
Then 
$$
\HP(\Omega^1_R)= 2\sum_{j=1}^s k_j -
\textstyle{\sum\limits_{\mathrm{char}(K)\nmid k_j}} 1.
$$
\end{corollary}

\begin{proof}
Let $\Y$ be the differential fattening of $\X$. 	
By Corollary~\ref{Cor-HFOmega^1_R}, we have 
	\begin{align*}
		\HP(\Omega^1_R) &=(n+2)\deg(\X) -\deg(\Y) \\
		&=(n+2)\dim_K(S) - \dim_K(A/J_\X^\diff)\\
		&= (n+2)\sum_{j=1}^s\dim_K(\mathcal{O}_j) -\sum_{j=1}^s\dim_K(A/J_{p_j}^\diff)\\
		&= \sum_{j=1}^s[(n+2)k_j -\dim_K(A/J_{p_j}^\diff)].
	\end{align*}
	Fix $j\in\{1,\dots,s\}$. Without loss of generality, 
	we may assume that $J_{p_j} = I_{p_j}^\deh = 
	\langle X_1,\dots,X_{n-1}, X_n^{k_j} \rangle$. 
	If $\charac(K) \nmid k_j$ then Corollary~\ref{Cor-Diff-MonomialIdeal}(b) implies that
	$$
	J_{p_j}^\diff = \langle X_1,\dots,X_{n-1}\rangle^2 + 
	\langle X_iX_n^{k_j} \mid i=1,\dots,n\rangle,
	$$
	and hence $\dim_K(A/J_{p_j}^\diff) = nk_j+1$.
	If $\charac(K) \mid k_j$ then Corollary~\ref{Cor-Diff-MonomialIdeal}(b) yields that
	$$
	J_{p_j}^\diff = \langle X_1,\dots,X_{n-1}\rangle^2 + 
	\langle X_n^{k_j} \rangle,
	$$
	and consequently $\dim_K(A/J_{p_j}^\diff) = nk_j$. Thus we obtain the above 
	formula for the constant Hilbert polynomial $\HP(\Omega^1_R)$.
\end{proof}


\bigbreak
%
%

\section{Differential Modules for 0-Dimensional LMG-Schemes}
\label{Section 4}

In the following, we continue with the established notation.
More precisely, let $K$ be a perfect field, and let $\X$ be a 0-dimensional scheme 
in $\mathbb{P}^n_K$ with support $\Supp(\X)=\{p_1,\dots,p_s\}$.
Its homogeneous vanishing ideal in $P=K[X_0,\dots,X_n]$ is $I_\X$
and its homogeneous coordinate ring is $R=P/I_\X$.
Furthermore, let $J_{p_j}$ be the affine vanishing ideal
of $p_j$ in $A=K[X_1,\dots,X_n]$ for $j=1,\dots,s$. 
The affine vanishing ideal of~$\X$ is $J_\X=J_{p_1}\cap\cdots\cap J_{p_s}$,
and its affine coordinate ring is $S=R/\langle x_0-1\rangle=A/J_\X$.
The ring $S$ can be decomposed as 
$S= \mathcal{O}_1\times\cdots\times\mathcal{O}_s$, where 
$\mathcal{O}_j= A/J_{p_j}$ 
is the 0-dimensional local ring of $\X$ at $p_j$.

In this section, we turn our attention to the following special class of 
0-dimen\-sional schemes that generalizes curvilinear schemes.

\begin{definition}
	\begin{enumerate}
		\item[(a)] The scheme $\X$ is called a \textbf{locally monomial scheme} if, 
		for each $j\in\{1,\dots,s\}$, there is a linear change of coordinates 
		$\phi_j: A \rightarrow A$ such that $\mathfrak{Q}_j=\phi_j(J_{p_j})$ is 
		a monomial ideal in $A$.
		
		\item[(b)] The scheme $\X$ is called a \textbf{locally monomial Gorenstein scheme}
		(or, for short, an \textbf{LMG-scheme}), if $\X$ is a locally monomial scheme
		and its local rings $\mathcal{O}_j \cong A/\mathfrak{Q}_j$ are Gorenstein
		for $j=1,\dots,s$.
	\end{enumerate}
\end{definition}

It is well-known that every locally complete intersection is locally Gorenstein;
however, the converse does not hold in general.
The following characterization of an LMG-scheme at a given point can be found in
\cite[Corollary 1.3.2 and Proposition A.6.5]{HH2011}.

\begin{proposition} 
	Let $\mathfrak{Q}$ be a 0-dimensional monomial ideal in $A\!=K[X_1,\dots,X_n]$.
	Then the following conditions are equivalent.
	\begin{enumerate}
		\item[(a)] The ring $A/\mathfrak{Q}$ is a Gorenstein ring.
		
		\item[(b)] The ring $A/\mathfrak{Q}$ is a complete intersection.
		
		\item[(c)] The ideal $\mathfrak{Q}$ is irreducible.
		
		\item[(d)] There are positive integers $k_1,\dots,k_n$ such that 
		$\mathfrak{Q} = \langle X_1^{k_1},\dots, X_n^{k_n}\rangle$. 
	\end{enumerate}
\end{proposition}

\begin{corollary}\label{Cor-Describing LMG-Schemes}
	The scheme $\X$ is an LMG-scheme if and only if, 
	for $j\in\{1,\dots,s\}$, there are positive integers $k_1,\dots,k_{n}$ 
	such that $\mathcal{O}_j \cong A/\langle X_1^{k_1},\dots, X_{n}^{k_n}\rangle$. 
\end{corollary}

Now we aim to establish a formula for the Hilbert polynomial of $\Omega^{m}_{R}$,
the module of K\"ahler differential $m$-forms associated to LMG-schemes $\X$. 
Since $\Omega^{m}_{R} =0$ for all $m>n+1$, we only need to consider 
the range $1\le m\le n+1$.
When $m=1$, we can apply Corollary~\ref{Cor-Diff-MonomialIdeal} to derive
the following formula for $\HP(\Omega^1_R)$.

\begin{proposition}\label{Prop-HPOfOmega^1_R-LMG-Scheme}
Let $\X$ be a 0-dimensional LMG-scheme in $\mathbb{P}^n_K$ with support
$\Supp(\X)=\{p_1,\dots,p_s\}$, and let $p$ denote the characteristic of~$K$.
For $j=1,\dots s$, we write $\mathcal{O}_j \cong A/\mathfrak{Q}_j$,
where $\mathfrak{Q}_j = \langle X_1^{k_{1j}},\dots,X_n^{k_{nj}} \rangle$ 
with $k_{ij}\ge 1$.
Then 
$$
\HP(\Omega^1_R) = 
\sum_{j=1}^s[(n+1)k_{1j}\cdots k_{nj} - \sum_{p\nmid k_{ij}}k_{1j}\cdots \widehat{k_{ij}}\cdots k_{nj}].
$$
\end{proposition}

\begin{proof}
Let $\Y$ be the differential fattening of $\X$. 	
By Corollary~\ref{Cor-HFOmega^1_R}, we have 
\begin{align*}
	\HP(\Omega^1_R) &=(n+2)\deg(\X) -\deg(\Y) \\
	&= (n+2)\sum_{j=1}^s\dim_K(\mathcal{O}_j) 
	-\sum_{j=1}^s\dim_K(A/\mathfrak{Q}_j^\diff)\\
	&= \sum_{j=1}^s[(n+2)k_{1j}\cdots k_{nj} -\dim_K(A/\mathfrak{Q}_j^\diff)].
\end{align*}
Letting $j\in\{1,\dots,s\}$, we will compute $\dim_K(A/\mathfrak{Q}_j^\diff)$.
Without loss of generality, we may assume that $p \mid k_{ij}$ for $i\le e$
and $p\nmid k_{ij}$ for $i>e$.
An application of Corollary~\ref{Cor-Diff-MonomialIdeal}(b) gives 
$$
\mathfrak{Q}_j^\diff = \langle X_1^{k_{1j}},\dots, X_e^{k_{ej}},
X_{e+1}^{k_{e+1 j}+1},\dots, X_n^{k_{n j}+1}\rangle
+ \langle X_{i_1}^{k_{i_1j}}X_j^{k_{i_2j}} \mid e< i_1<i_2\le n\rangle.
$$
Consider a term $t=X_1^{\alpha_1}\cdots X_n^{\alpha_n} \in \mathbb{T}^{n+1}$.
Note that $t\notin \mathfrak{Q}_j^\diff$ yields 
$\alpha_i \ne k_{ij}$ for $i=1,\dots,e$.
Thus, $t\in \mathbb{T}^{n+1} \setminus \mathfrak{Q}_j^\diff$ if and only if
one of the following holds:
\begin{itemize}
	\item $0\le \alpha_i < k_{ij}$ for all $i=1,\dots,n$; 
	\item if $\alpha_i = k_{ij}$ for some $i\in\{e+1,\dots,n\}$, 
	then $0\le \alpha_l < k_{lj}$ for all $l\ne i$. 
\end{itemize}
It follows that 
\begin{align*}
\dim_K(A/\mathfrak{Q}_j^\diff) \;&=\; \#(\mathbb{T}^{n+1}\setminus \mathfrak{Q}_j^\diff)
= k_{1j}\cdots k_{nj}+\sum_{i=e+1}^nk_{1j}\cdots \widehat{k_{ij}}\cdots k_{nj}\\
&=k_{1j}\cdots k_{nj}+\sum_{p\nmid k_{ij}}k_{1j}\cdots \widehat{k_{ij}}\cdots k_{nj}.
\end{align*}
where $\widehat{k_{ij}}$ means that the term $k_{ij}$ is omitted from the product. 
Consequently, we get the formula
\begin{align*}
\HP(\Omega^1_R) &= \sum_{j=1}^s[(n+2)k_{1j}\cdots k_{nj} - 
\dim_K(A/\mathfrak{Q}_j^\diff)]\\
&= \sum_{j=1}^s[(n+1)k_{1j}\cdots k_{nj} - \sum_{p\nmid k_{ij}}k_{1j}\cdots \widehat{k_{ij}}\cdots k_{nj}].
\end{align*}
\end{proof}

\begin{example}\label{Exam-ComputingHFOmega^1_R-LMG-Sch}
Let $K=\mathbb{F}_3$ and let $\X$ be a 0-dimensional scheme in $\mathbb{P}^3_K$ 
defined by $I_\X=\langle f_1,\dots,f_6\rangle \subseteq P=K[X_0,X_1,X_2,x_3]$, where 
\begin{align*}
f_1 &= X_{1}^3  +X_{2}^3,\\
f_2 &= X_{0}^2 X_{1} +X_{0} X_{1}^2  +X_{0}^2 X_{2} +X_{0} X_{1} X_{2} +X_{1}^2 X_{2} -X_{2}^3,\\
f_3 &= X_{0}^2 X_{2} X_{3} +X_{0} X_{1} X_{2} X_{3} +X_{1}^2 X_{2} X_{3} -X_{3}^4,\\
f_4 &= X_{0}^3 X_{3} +X_{2}^3 X_{3} +X_{3}^4,\\
f_5 &= X_{0}^3 X_{2} +X_{0}^2 X_{2}^2  +X_{0} X_{1} X_{2}^2  +X_{1}^2 X_{2}^2  +X_{2}^4,\\
f_6 &= X_{0}^2 X_{2}^3  +X_{0} X_{1} X_{2}^3  +X_{1}^2 X_{2}^3  -X_{0}^2 X_{3}^3  -X_{0} X_{1} X_{3}^3  -X_{1}^2 X_{3}^3.
\end{align*}
Then $\Supp(\X)=\{p_1,p_2,p_3\}$. In $A=K[X_1,X_2,X_3]$, 
the affine vanishing ideal $J_\X$ decomposes 
$J_\X = J_{p_1} \cap J_{p_2} \cap J_{p_3}$,
where $J_{p_1} = \langle X_1, X_2, X_3 \rangle$,
$J_{p_2} = \langle X_1+1, X_2-1, (X_3-1)^3 \rangle$, and
$J_{p_3} = \langle X_1^2+X_1+1, X_2^3+1, X_3^4 \rangle$.
Letting $\mathfrak{Q}_1 = \langle X_1, X_2, X_3 \rangle$, $\mathfrak{Q}_2
=\langle X_1, X_2, X_3^3 \rangle$, and 
$\mathfrak{Q}_3 =\langle X_1^2, X_2^3, X_3^4 \rangle$, we find 
$\mathcal{O}_j \cong A/\mathfrak{Q}_j$ for $j=1,2,3$.
Thus $\X$ is a 0-dimensional LMG-scheme.  
An application of Proposition~\ref{Prop-HPOfOmega^1_R-LMG-Scheme} yields
\begin{align*}
	\HP(\Omega^1_R) &= [(3+1) - 3] + [(3+1)3 -3-3] + [(3+1)\cdot 2\cdot 3\cdot 4 -
	3\cdot 2 - 3\cdot 4] \\
	&= 1 + 6 + 78 = 85.
\end{align*}
\end{example}

In order to treat the case that $2\le m\le n+1$, we apply the following observation.
By \cite[Prop. 4.12]{Kun1986}, the module $\Omega^{m}_{S}$ 
decomposes as 
$\Omega^{m}_{S} \;\cong\; \Omega^{m}_{\mathcal{O}_1}\times
\cdots\times \Omega^{m}_{\mathcal{O}_s}$.
Furthermore, according to Proposition~\ref{prop:HigherKDM},
the Hilbert polynomial of $\Omega^{m}_R$ satisfies 
\begin{align*}
	\HP(\Omega^{m}_R) &\;=\; \dim_K(\Omega^{m}_{S}) + \dim_K(\Omega^{m-1}_{S})\\
	&\;=\; \sum_{j=1}^s\; [\, \dim_K(\Omega^{m}_{\mathcal{O}_j}) 
	+ \dim_K(\Omega^{m-1}_{\mathcal{O}_j})\,].
\end{align*}	
and the regularity index of $\Omega^{m}_R$ is bounded by 
$\mathrm{ri}(\Omega^{m}_R) \le 2r_\X +m$.
This reduces the problem to computing $\dim_K(\Omega^{m}_{S})$ 
in the local setting where $\X$ contains exactly one point.

\begin{proposition}\label{Prop-DimOfOmegaSm}
Let $m,k_1,\dots,k_n$ be positive integers and $p=\charac(K)$.
Suppose that $S=A/\mathfrak{Q}$, where $A=K[X_1,\dots,X_n]$ and
$\mathfrak{Q} = \langle X_1^{k_{1}},\dots,X_n^{k_{n}} \rangle\subseteq A$.
For $1\le i_1<\cdots<i_m\le n$, let 
$\Gamma_{(i_1,\dots,i_m)} = \{i\in \{i_1,\dots,i_m\} :\;\,	p \nmid k_i\}$ 
and define 
$$
J_{(i_1,\dots,i_m)} = \mathfrak{Q} +
\langle X_{i}^{k_{i}-1} \mid i\in \Gamma_{(i_1,\dots,i_m)}\rangle.
$$
Then 
$$
\Omega^m_S \cong \bigoplus_{1\le i_1<\dots<i_m\le n} (A/J_{(i_1,\dots,i_m)})(-m)
$$
and
$$
\dim_K(\Omega^m_S) = \sum_{1\le i_1<\dots<i_m\le n} 
	\prod_{i\notin \Gamma_{(i_1,\dots,i_m)}}k_i 
	\prod_{i\in\Gamma_{(i_1,\dots,i_m)}} (k_i-1).
$$
\end{proposition}

\begin{proof}
We equip $A$ with the standard grading, where $\deg(X_i) = \deg(dX_i)=1$
for every $i=1,\dots,n$. Then $\Omega^m_A$ is a graded free $A$-module with basis 
$$
\{dX_{i_1}\wedge\cdots\wedge dX_{i_m} \mid 1\le i_1<\cdots<i_m\le n\}
$$
by Proposition~\ref{prop:PresentationOmegaRm}(b).
Let us set $e_{i_1,\dots,i_m} = dX_{i_1}\wedge\cdots\wedge dX_{i_m}$ and write 
$$
\Omega^m_A = \bigoplus_{1\le i_1<\dots<i_m\le n} Ae_{i_1,\dots,i_m}.
$$
By \cite[Proposition 4.12]{Kun1986}, we have the exact sequence of graded $A$-modules 
$$
0\To \mathfrak{Q}\Omega^m_A +d\mathfrak{Q}\wedge\Omega^{m-1}_A
\To \Omega^m_A \To \Omega^m_S \To 0.
$$
Since $\mathfrak{Q}$ is a monomial ideal, the graded $A$-module 
$\mathfrak{Q}\Omega^m_A +d\mathfrak{Q}\wedge\Omega^{m-1}_A$ is a monomial module 
in the sense of \cite[Definition 1.3.7]{KR2000}. Of course, for $i=1,\dots,n$,
we have $dX_i^{k_i} =\sum_{i=1}^n\partial_i(X_i^{k_i}) =k_i X_i^{k_i-1}dX_i$ 
if $p\nmid k_i$ and $dX_i^{k_i}=0$ if $p\mid k_i$.
By \cite[Theorem 1.3.9]{KR2000} and Proposition~\ref{prop:PresentationOmegaRm}(d), 
a calculation gives 
$$
\mathfrak{Q}\Omega^m_A +d\mathfrak{Q}\wedge\Omega^{m-1}_A 
= \bigoplus_{1\le i_1<\dots<i_m\le n} J_{(i_1,\dots,i_m)}e_{i_1,\dots,i_m}.
$$
It follows that
$$
\Omega^m_S= \bigoplus_{1\le i_1<\dots<i_m\le n} (A/J_{(i_1,\dots,i_m)})(-m).
$$
Moreover, we have $\dim_K(A/J_{(i_1,\dots,i_m)})
= \prod_{i\notin \Gamma_{(i_1,\dots,i_m)}}k_i 
\prod_{i\in\Gamma_{(i_1,\dots,i_m)}} (k_i-1)$.
Consequently, we obtain the above formula for 
$\dim_K(\Omega^m_S)$, as desired.
\end{proof}

\begin{example} 
Let us continue looking at the LMG-scheme $\X$ given 
in Example~\ref{Exam-ComputingHFOmega^1_R-LMG-Sch}.
We know that $\Supp(\X)=\{p_1,p_2,p_3\}$ and 
$\mathcal{O}_j \cong A/\mathfrak{Q}_j$ for $j=1,2,3$,
where $\mathfrak{Q}_1 = \langle X_1, X_2, X_3 \rangle$, $\mathfrak{Q}_2
=\langle X_1, X_2, X_3^3 \rangle$, and 
$\mathfrak{Q}_3 =\langle X_1^2, X_2^3, X_3^4 \rangle$.
Since $K$ is a perfect field and $\mathcal{O}_1\cong K$,
we have $\Omega^m_{\mathcal{O}_1} =0$ for $m\ge 1$.
Notice that $\Omega^m_{\mathcal{O}_j}=0$ for $m\ge 4$ and $j=2,3$.
By Proposition~\ref{Prop-DimOfOmegaSm}, we find
\begin{align*}
\dim_K(\Omega^1_{\mathcal{O}_2}) &= 3, 
&&	\dim_K(\Omega^1_{\mathcal{O}_3}) = 54,\\
\dim_K(\Omega^1_{\mathcal{O}_2}) &= 0,
&& \dim_K(\Omega^2_{\mathcal{O}_3}) = 39,\\
\dim_K(\Omega^1_{\mathcal{O}_2}) &= 0, 
&&	\dim_K(\Omega^3_{\mathcal{O}_3}) = 9.
\end{align*}
Consequently, by Propositions~\ref{prop:HigherKDM}, we get
\begin{align*}
	\HP(\Omega^1_R) &=\deg(\X) + \sum_{j=1}^{3}\dim_K\Omega^1_{\mathcal{O}_j}
	= 28 + 0+3+54 = 85, \\
	\HP(\Omega^2_R) &=\sum_{j=1}^{3}(\dim_K\Omega^1_{\mathcal{O}_j}+\dim_K\Omega^2_{\mathcal{O}_j})
	= 0 +3 +93 = 96, \\
	\HP(\Omega^3_R) &=\sum_{j=1}^{3}(\dim_K\Omega^2_{\mathcal{O}_j}+\dim_K\Omega^3_{\mathcal{O}_j})
	= 0 + 0 + 48 = 48.
\end{align*}
\end{example}

When $\X$ is a 0-dimensional LMG-scheme, an application 
of Propositions~\ref{prop:HigherKDM} and~\ref{Prop-DimOfOmegaSm} yields 
the following formula for the Hilbert polynomial of $\Omega^m_R$.

\begin{corollary}\label{Cor-HPofOmage^m_R-LMG-Sch}
Let $\X$ be a 0-dimensional LMG-scheme in $\mathbb{P}^n_K$, and let $p=\mathrm{char}(K)$.
For $j=1,\dots s$, we write $\mathcal{O}_j \cong A/\mathfrak{Q}_j$,
where $\mathfrak{Q}_j = \langle X_1^{k_{1j}},\dots,X_n^{k_{nj}} \rangle$ 
with $k_{ij}\ge 1$. For $1\le i_1<\cdots<i_m\le n$ and $1\le j\le s$, let 
$$
\Gamma_{(i_1,\dots,i_m),j} = \{i\in \{i_1,\dots,i_m\} :\;\,	p \nmid k_{ij}\}.
$$
For $m\ge 1$, we have 
\begin{align*}
\HP(\Omega^m_R) &= \sum_{j=1}^s
\Big[ \sum_{1\le i_1<\dots<i_{m-1}\le n} 
\prod_{i\notin \Gamma_{(i_1,\dots,i_{m-1}),j}}k_{ij} 
\prod_{i\in\Gamma_{(i_1,\dots,i_{m-1}),j}} (k_{ij}-1) \\
&\quad +\sum_{1\le l_1<\dots<l_m\le n} 
\prod_{i\notin \Gamma_{(l_1,\dots,l_m),j}}k_{ij}
\prod_{i\in\Gamma_{(l_1,\dots,l_m),j}} (k_{ij}-1)\Big].
\end{align*}
\end{corollary}

%
%

\section*{Acknowledgements}

The authors would like to thank Martin Kreuzer for his generous support 
and encouragement. This work is supported by 
the \textit{Vietnam Ministry of Education and Training} 
under the grant number B2025-DHH-02. 

%
%

\end{document}